# Timelike Minimal Surfaces via Loop Groups

*The first named author would like to dedicate this paper to the memory of his father Yukio*

J. INOGUCHI[1],[★] and M. TODA[2]

[1]*Department of Mathematics Education, Faculty of Education, Utsunomiya University, Utsunomiya, 321-8505, Japan. e-mail: inoguchi@cc.utsunomiya-u.ac.jp*

[2]*Department of Mathematics and Statistics, Texas Tech University, Lubbock, TX 79407, U.S.A. e-mail: mtoda@ttu.edu*

**Abstract.** This work consists of two parts. In Part I, we shall give a systematic study of Lorentz conformal structure from structural viewpoints. We study manifolds with split-complex structure. We apply general results on split-complex structure for the study of Lorentz surfaces.

In Part II, we study the conformal realization of Lorentz surfaces in the Minkowski 3-space via conformal minimal immersions. We apply loop group theoretic Weierstrass-type representation of timelike constant mean curvature for timelike minimal surfaces. Classical integral representation formula for timelike minimal surfaces will be recovered from loop group theoretic viewpoint.

**Mathematics Subject Classifications (2000):** 53A10, 58E20, 22E67.

**Key words:** Lorentz surfaces, harmonic maps, loop groups, timelike minimal surfaces.

**Introduction**

The classical wave equation:

$$-u_{tt} + u_{xx} = 0$$

over a plane $\mathbb{R}^2(t, x)$ can be understood as the *harmonicity equation* with respect to the Lorentz metric $-dt^2 + dx^2$. It is easy to see that Lorenz harmonicity is invariant under conformal transformation of $(\mathbb{R}^2(t, x), -dt^2 + dx^2)$.

This observation tells us that for the study of second order PDE of hyperbolic type, it is convenient to introduce the notion of *Lorentz conformal structure*.

Oriented 2-manifolds together with Lorenz conformal structure are called *Lorentz surfaces*.

R. Kulkarni [42] initiated the global study of Lorentz conformal structure and Lorentz surfaces. (For former studies on local structure, see Burns [7], Dzan [20], and Smith [59].) T. Weinstein (also known as T. K. Milnor) and her collaborators

---

[★] Partially supported by MEXT. Grant-in-Aid for Young Scientists 14740053.

further analyzed Lorentz surfaces [38, 45, 46, 51–53, 60–63, 66, 67]. In particular, they obtained many interesting results on the global behavior of boundaries of Lorentz surfaces. In particular, Smyth [60, 61] studied completions of conformal boundaries. Klarreich [38] studied smoothability of conformal boundaries.

The authors mentioned above investigated how Lorentz surfaces can be realized as *timelike minimal surfaces* in Minkowski 3-space.

On the other hand, Luo and Stong studied conformal imbeddings of discs with Lorentzian metric into Minkowski plane [47].

Weinstein showed that Bernstein problems for timelike minimal surfaces should be considered in conformal setting. (For a precise setting see [51].) The *conformal Bernstein problem* in $\mathbb{E}^3_1$ was solved independently by Weinstein [51] and Magid [49]. (See also further generalizations, [46].) McNertney gave a Weierstrass formula for timelike minimal surfaces in her (unpublished) doctoral dissertation [50]. She classified timelike minimal surfaces of revolution and their *associated helicoids*. Later Van de Woestijine [65] independently gave a similar classification table of timelike minimal surfaces of revolution. Moreover he classified timelike minimal ruled surfaces and timelike minimal translation surfaces.

Gu [31], Magid [48] and Weinstein [67] gave several variants of Weierstrass representation formulas. Gu [31] discussed relations between timelike minimal surface theory and fluid dynamics.

The study of timelike minimal surfaces has yet another motivation. The fundamentals of classical string theory can be summarized as follows: A *closed string* is an object $\gamma$ in the physical space, that is homeomorphic to $S^1$. Intuitively speaking, a string evolves in time while sweeping a surface $\Sigma$, called *world sheet*, in spacetime. For physical reasons $\Sigma$ is supposed to be a timelike surface. The dynamical equations for a string are defined by a variational principle: The first area variation of $\Sigma$ must vanish subject to the condition that the initial and final configuration of the string are kept fixed. Hence $\Sigma$ will be a timelike minimal surface having two spacelike boundary components. For more details, we refer to [9, 41].

Cauchy problems for timelike minimal surfaces in spacetimes are studied by Deck [13, 14] and Gu [30].

It is easy to see that timelike minimal surfaces in a spacetime can be realized via conformal harmonic maps of Lorentz surfaces. Conformal harmonic maps of a Lorentz surface into a semi-Riemannian symmetric space $G/K$ are called *nonlinear sigma models* (with symmetry group $G$) in particle physics. Nonlinear sigma models are regarded as toy models of gauge theory. Harmonic maps of Lorentz surfaces into semi-Riemannian manifolds are often called *wave maps*.

Harmonic maps of Lorentz surfaces into the 2-sphere $S^2$ are closely related to soliton theory. In fact, such harmonic maps are viewed as Gauss maps of pseudospherical surfaces in Euclidean 3-space. More precisely, the second fundamental form of a pseudospherical surface determines a Lorentz conformal structure on the surface (frequently called the *second conformal structure*). With respect to the second conformal structure, the Gauss map of the pseudospherical surface is a harmonic map. Conversely, for every given (weakly regular) harmonic map of a Lorentz surface, there exists a pseudospherical surface, having this harmonic map as a Gauss map. Such surfaces can be given explicitly by *Lelieuvre's formula* [43].

In our previous work [34], the first named author established a $2 \times 2$ matrix setting for timelike surfaces with constant mean curvature.

Dorfmeister, Pedit and the second named author investigated minimal surfaces via a loop group method (a special case of the method frequently called *DPW-method*) [18].

Recently, Dorfmeister and the present authors gave loop group theoretic Weierstraß-type representation (or *nonlinear d'Alemberet formula*) for timelike constant mean curvature surfaces [17]. In this paper we use loop group techniques in order to construct timelike minimal surfaces.

The present work consists of two parts. In Part I, we shall give a systematic study of Lorentz conformal structures from a structural viewpoint. We start with split-complex numbers and split-quaternions (Section 1). These ingredients will be used for the study of a conformal realization of Lorentz surfaces in Minkowski 3-space.

In Section 2, we study manifolds with split-complex structure, especially split Kähler structure. In Section 3, we shall apply the general results of Sections 1–2, to Lorentz surfaces.

In Part II, we study the conformal realization of Lorentz surfaces in the Minkowski 3-space. Our main interest lies in minimal ones. We start with a loop group approach, similar to the one used for timelike constant mean curvature surfaces developed in [17]. We show that our scheme in [17] can be also applied to timelike minimal surfaces. Then we deduce the classical Weierstrass formula for the minimal case.

## Part I. Lorentz Conformal Structure

## 1. Split-Complex Numbers and Split-Quaternions

### 1.1. SPLIT-COMPLEX NUMBERS

Let $\mathbb{C}'$ be a real algebra spanned by 1 and $k'$, with multiplication law:

$$1 \cdot k' = k' \cdot 1 = k', \quad (k')^2 = 1. \tag{1.1}$$

The algebra $\mathbb{C}' = \mathbb{R}1 \oplus \mathbb{R}k'$ is commutative. An element of $\mathbb{C}'$ is called a *split-complex number*. The algebra $\mathbb{C}'$ is called the *algebra of split-complex numbers*. Note that split-complex numbers are also called *paracomplex numbers* (or *hyperbolic numbers* [6]). Let $\zeta$ be a split-complex number. Then $\zeta$ can be uniquely expressed as $\zeta = x + yk'$. The *real part* $\operatorname{Re} \zeta$ and the *imaginary part* $\operatorname{Im} \zeta$ of $\zeta$ are defined by

$$\operatorname{Re} \zeta := x, \quad \operatorname{Im} \zeta := y.$$

For a split-complex number $\zeta = x + yk'$, the *conjugate* $\bar{\zeta}$ is defined by $\bar{\zeta} := x - yk'$. Since $-\zeta\bar{\zeta} = -x^2 + y^2$, the algebra $\mathbb{C}'$ is identified with the Minkowski plane $\mathbb{E}_1^2$:

$$\mathbb{E}_1^2(x, y) = (\mathbb{R}^2(x, y), -dx^2 + dy^2)$$

as a scalar product space. Namely we identify $\mathbb{C}'$ with $\mathbb{E}_1^2$ via

$$1 \longleftrightarrow \begin{pmatrix} 1 \\ 0 \end{pmatrix}, \qquad k' \longleftrightarrow \begin{pmatrix} 0 \\ 1 \end{pmatrix}.$$

The scalar product $\langle \cdot, \cdot \rangle$ of $\mathbb{E}_1^2$ can be written as

$$\langle \zeta_1, \zeta_2 \rangle = -\mathrm{Re}(\bar{\zeta}_1 \zeta_2)$$

for $\zeta_1 = x_1 + y_1 k'$, $\zeta_2 = x_2 + y_2 k' \in \mathbb{C}'$.

The "unit" circle group $\mathrm{U}(1; \mathbb{C}')$ of $\mathbb{C}'$:

$$\mathrm{U}(1; \mathbb{C}') = \{\zeta \in \mathbb{C}' \mid \zeta\bar{\zeta} = 1\}$$

is the hyperbola

$$H_0^1 = \{(x, y) \in \mathbb{E}_1^2 \mid -x^2 + y^2 = -1\}$$

which is isomorphic to the multiplicative group $\mathbb{R}^\times = (\mathbb{R} \setminus \{0\}, \times)$.

## 1.2. SPLIT COMPLEX STRUCTURE

The multiplication by $k'$ on $\mathbb{C}'$ defines a linear endomorphism $\mathcal{K}'$:

$$\mathcal{K}'1 = k', \qquad \mathcal{K}'k' = 1. \tag{1.2}$$

Clearly, this linear endomorphism satisfies

$$\mathcal{K}' \mathcal{K}' = \text{identity}, \qquad \langle \mathcal{K}'\zeta_1, \mathcal{K}'\zeta_2 \rangle = -\langle \zeta_1, \zeta_2 \rangle.$$

The linear anti isometry $\mathcal{K}'$ is called the *associated split-complex structure* of $\mathbb{C}'$.

## 1.3. SPLIT-COMPLEX LINEAR SPACES

More generally, we shall consider split-complex structures on linear spaces.

DEFINITION 1.1. Let $V$ be a finite-dimensional real linear space $V$. A linear endomorphism $\mathcal{F}$ on $V$ is said to be a *product structure* on $V$ if $\mathcal{F}^2 = $ identity.

Since the eigenvalues of $\mathcal{F}$ are $\pm 1$, $V$ can be decomposed as

$$V = V_+ \oplus V_-.$$

Here $V_+$ and $V_-$ are the eigenspaces corresponding to $1$ and $-1$ respectively.

Let $\mathcal{K}'$ be a product structure such that $\dim V_+ = \dim V_- \geq 1$. Then $\mathcal{K}'$ defines a left action of $\mathbb{C}'$ on $V$:

$$(a + bk')\mathbf{x} := a\mathbf{x} + b\mathcal{K}'\mathbf{x}.$$

Thus $(V, \mathcal{K}')$ is regarded as a free $\mathbb{C}'$-module. This fact motivates the following definition:

DEFINITION 1.2. Let $V$ be a finite-dimensional real linear space. A product structure $\mathcal{K}'$ on $V$ is said to be a *split-complex structure* if

$$\dim V_+ = \dim V_-.$$

In this paper we restrict our attention to split-complex structures.

*Remark 1.3.* Some authors use the terminology *paracomplex* instead of *split-complex* ([11, 21, 22, 35, 36]).

A vector in a split-complex linear space $(V, \mathcal{K}')$ is said to be a vector of $(1, 0)$-*type* [resp. $(0, 1)$-*type*] if it is in $V_+$ [resp. $V_-$].

The split-complex structure $\mathcal{K}'$ naturally induces a split-complex structure on the dual linear space $V^*$ of $V$:

$$(\mathcal{K}'A)(\mathbf{x}) := A(\mathcal{K}'\mathbf{x}), \quad \mathbf{x} \in V, \ A \in V^*.$$

Denote by $V'$, $V''$ the $\pm 1$-eigenspaces of $\mathcal{K}'$ on $V^*$. Then $V^*$ splits into:

$$V^* = V' \oplus V''.$$

A covector in $V^*$ is said to be *of $(1, 0)$-type* [resp. $(0, 1)$-type] if it is an element of $V'$ [resp. $V''$].

PROPOSITION 1.4. *Let $(V_1, \mathcal{K}'_1)$, $(V_2, \mathcal{K}'_2)$ be split-complex linear spaces and $F: V_1 \to V_2$ a linear map. Then the following conditions are mutually equivalent*:

(1) *$F$ is $\mathbb{C}'$-linear*;
(2) *$F \circ \mathcal{K}'_1 = \mathcal{K}'_2 \circ F$*;
(3) *$F((V_1)_+) \subset (V_2)_+$, $F((V_1)_-) \subset (V_2)_-$*;
(4) *$F^*(V'_2) \subset V'_1$, $F^*(V''_2) \subset V''_1$, where $F^*$ denotes dual mapping of $F$.*

COROLLARY 1.5. *Let $(V_1, \mathcal{K}'_1)$, $(V_2, \mathcal{K}'_2)$ be split-complex linear spaces and $F: V_1 \to V_2$ a linear map. Then the following conditions are mutually equivalent*:

(1) *$F$ is anti $\mathbb{C}'$-linear*;
(2) *$F \circ \mathcal{K}'_1 = -\mathcal{K}'_2 \circ F$*;

(3) $F((V_1)_+) \subset (V_2)_-$, $F((V_1)_-) \subset (V_2)_+$;
(4) $F^*(V_2') \subset V_1''$, $F^*(V_2'') \subset V_1'$.

DEFINITION 1.6. Let $(V, \mathcal{K}')$ be a split-complex linear space. A scalar product $\langle \cdot, \cdot \rangle$ on $V$ is said to be a *split-Hermitian scalar product* (or *para-Hermitian scalar product*) if

$$\langle \mathcal{K}'\mathbf{x}, \mathcal{K}'\mathbf{y} \rangle = -\langle \mathbf{x}, \mathbf{y} \rangle$$

for all $\mathbf{x}, \mathbf{y} \in V$. A triple $(V, \mathcal{K}', \langle \cdot, \cdot \rangle)$ is called a *split-Hermitian linear space* or *para-Hermitian linear space*.

Obviously, the algebra $\mathbb{C}'$ of split-complex numbers is a typical example of split-Hermitian linear space.

### 1.4. SPLIT-QUATERNIONS

In the end of this section, we give an example of real 4-dimensional split-Hermitian linear space.

Let us denote the Clifford algebra of $\mathbb{C}'$ by $\mathbb{H}'$. The Clifford algebra $\mathbb{H}'$ is isomorphic to the following real algebra:

$$\mathbb{R}\mathbf{1} \oplus \mathbb{R}\mathbf{i} \oplus \mathbb{R}\mathbf{j}' \oplus \mathbb{R}\mathbf{k}'$$

generated by the basis $\{\mathbf{1}, \mathbf{i}, \mathbf{j}', \mathbf{k}'\}$;

$$\mathbf{1}\mathbf{i} = \mathbf{i}\mathbf{1} = \mathbf{i}, \quad \mathbf{1}\mathbf{j}' = \mathbf{j}'\mathbf{1} = \mathbf{j}', \quad \mathbf{1}\mathbf{k}' = \mathbf{k}'\mathbf{1} = \mathbf{k}',$$

$$\mathbf{i}\mathbf{j}' = -\mathbf{j}'\mathbf{i} = \mathbf{k}', \quad \mathbf{j}'\mathbf{k}' = -\mathbf{k}'\mathbf{j}' = -\mathbf{i}, \quad \mathbf{k}'\mathbf{i} = -\mathbf{i}\mathbf{k}' = \mathbf{j}',$$

$$\mathbf{i}^2 = -\mathbf{1}, \quad \mathbf{j}'^2 = \mathbf{k}'^2 = \mathbf{1}.$$

An element of the Clifford algebra $\mathbb{H}'$ is called a *split-quaternion*. The Clifford algebra $\mathbb{H}'$ is called the algebra of split-quaternions (cf. [32, 57, 68, 69]).

We identify $\mathbb{C}'$ with $\mathbb{R}\mathbf{1} \oplus \mathbb{R}\mathbf{k}'$:

$$x1 + yk' \longmapsto x\mathbf{1} + y\mathbf{k}'.$$

The multiplication by $\mathbf{k}'$ defines a split-complex structure $\mathcal{K}'$ on $\mathbb{H}'$:

$$\mathcal{K}'\xi = \mathbf{k}'\xi, \quad \xi \in \mathbb{H}'.$$

Similarly, $\mathbf{j}'$ defines another split-complex structure $\mathcal{J}'$ on $\mathbb{H}'$:

$$\mathcal{J}'\xi = \mathbf{j}'\xi, \quad \xi \in \mathbb{H}'.$$

The algebra of split-quaternions $\mathbb{H}'$ is a 2-dimensional $\mathbb{C}'$-module:

$$\mathbb{H}' = \mathbb{C}' \oplus \mathbf{i}\mathbb{C}'.$$

Note that for any split-complex number $\zeta$,

$$\zeta \mathbf{i} = \mathbf{i}\bar{\zeta}.$$

On the other hand, the multiplication by $\mathbf{i}$ defines a complex structure $\mathcal{J}$ on $\mathbb{H}'$:

$$\mathcal{J}\xi = \mathbf{i}\xi, \quad \xi \in \mathbb{H}'.$$

Thus, we may identify the field of all complex numbers $\mathbb{C}$ with $\mathbb{R}\mathbf{1} \oplus \mathbb{R}\mathbf{i}$.

$$x1 + y\sqrt{-1} \longmapsto x\mathbf{1} + y\mathbf{i}.$$

Under this identification, $\mathbb{H}'$ is regarded as a right complex linear space:

$$\mathbb{H}' = \mathbb{C} \oplus \mathbf{k}'\mathbb{C}.$$

Note that for any complex number $z \in \mathbb{C}$,

$$z\mathbf{k}' = \mathbf{k}'\bar{z}.$$

1.5. SCALAR PRODUCT ON SPLIT-QUATERNION ALGEBRA

For a split-quaternion $\xi = \xi_0 \mathbf{1} + \xi_1 \mathbf{i} + \xi_2 \mathbf{j}' + \xi_3 \mathbf{k}'$, the *conjugate* $\bar{\xi}$ is defined by

$$\bar{\xi} = \xi_0 \mathbf{1} - \xi_1 \mathbf{i} - \xi_2 \mathbf{j}' - \xi_3 \mathbf{k}'.$$

It is easy to see that $-\xi\bar{\xi} = -\xi_0^2 - \xi_1^2 + \xi_2^2 + \xi_3^2$. Hereafter we identify $\mathbb{H}'$ with a semi-Euclidean space $\mathbb{E}_2^4$:

$$\mathbb{E}_2^4 = (\mathbb{R}^4(\xi_0, \xi_1, \xi_2, \xi_3), \ -d\xi_0^2 - d\xi_1^2 + d\xi_2^2 + d\xi_3^2).$$

Then it is easy to see that $(\mathbb{H}', \mathcal{K}', \langle \cdot, \cdot \rangle)$ is a 4-dimensional split-Hermitian linear space.

Let $G = \{\xi \in \mathbb{H}' \mid \xi\bar{\xi} = 1\}$ be the multiplicative group of timelike unit split-quaternions. The Lie algebra $\mathfrak{g}$ of $G$ is the imaginary part of $\mathbb{H}'$, that is,

$$\mathfrak{g} = \operatorname{Im} \mathbb{H}' = \{\xi_1 \mathbf{i} + \xi_2 \mathbf{j}' + \xi_3 \mathbf{k}' \mid \xi_1, \xi_2, \xi_3 \in \mathbb{R}\}.$$

The Lie bracket of $\mathfrak{g}$ is simply the commutator of the split-quaternion product. Note that the commutation relations of $\mathfrak{g}$ are given by:

$$[\mathbf{i}, \mathbf{j}'] = 2\mathbf{k}', \qquad [\mathbf{j}', \mathbf{k}'] = -2\mathbf{i}, \qquad [\mathbf{k}', \mathbf{i}] = 2\mathbf{j}'.$$

The Lie algebra $\mathfrak{g}$ is naturally identified with the Minkowski 3-space

$$\mathbb{E}_1^3 = (\mathbb{R}^3(\xi_1, \xi_2, \xi_3), \ -d\xi_1^2 + d\xi_2^2 + d\xi_3^2)$$

as a metric linear space.

## 1.6. MATRIX GROUP MODEL OF $\mathbb{H}'$

Next, we introduce a $2 \times 2$ matricial expression of $\mathbb{H}'$ as follows:

$$\xi = \xi_0 \mathbf{1} + \xi_1 \mathbf{i} + \xi_2 \mathbf{j}' + \xi_3 \mathbf{k}' \longleftrightarrow \begin{pmatrix} \xi_0 - \xi_3 & -\xi_1 + \xi_2 \\ \xi_1 + \xi_2 & \xi_0 + \xi_3 \end{pmatrix}. \quad (1.3)$$

In particular, the matricial expressions of the natural basis of $\mathbb{H}'$ are given by

$$\mathbf{1} \longleftrightarrow \begin{pmatrix} 1 & 0 \\ 0 & 1 \end{pmatrix}, \qquad \mathbf{i} \longleftrightarrow \begin{pmatrix} 0 & -1 \\ 1 & 0 \end{pmatrix},$$

$$\mathbf{j}' \longleftrightarrow \begin{pmatrix} 0 & 1 \\ 1 & 0 \end{pmatrix}, \qquad \mathbf{k}' \longleftrightarrow \begin{pmatrix} -1 & 0 \\ 0 & 1 \end{pmatrix}.$$

This correspondence gives an algebra isomorphism between $\mathbb{H}'$ and the algebra $M_2\mathbb{R}$ of all real $2 \times 2$ matrices. Under the identification (1.3), the group $G$ of timelike unit split-quaternions corresponds to the special linear group:

$$\mathrm{SL}_2\mathbb{R} = \left\{ \begin{pmatrix} a & b \\ c & d \end{pmatrix} \in M_2\mathbb{R} \mid ad - bc = 1 \right\}.$$

The semi-Euclidean metric of $\mathbb{H}'$ corresponds to the following scalar product on $M_2\mathbb{R}$:

$$\langle X, Y \rangle = \frac{1}{2}\{\mathrm{tr}(XY) - \mathrm{tr}(X)\mathrm{tr}(Y)\} \quad (1.4)$$

for all $X, Y \in M_2\mathbb{R}$. The metric of $G$ induced by (1.4) is a bi-invariant Lorentz metric of constant curvature $-1$. Hence the Lie group $G$ is identified with an anti de-Sitter 3-space $H_1^3$ of constant curvature $-1$. (See [55] and Dajzcer and Nomizu [12].)

As in the case of Euclidean 3-space, the *vector product operation* $\times$ of $\mathbb{E}_1^3$ is defined by

$$\xi \times \eta = (\xi_3\eta_2 - \xi_2\eta_3,\ \xi_3\eta_1 - \xi_1\eta_3,\ \xi_1\eta_2 - \xi_2\eta_1)$$

for $\xi = (\xi_1, \xi_2, \xi_3)$, $\eta = (\eta_1, \eta_2, \eta_3) \in \mathbb{E}_1^3$. If we regard $\xi$ and $\eta$ as elements of $\mathfrak{g}$, then the vector product of $\xi$ and $\eta$ is written in terms of the Lie bracket as follows:

$$\xi \times \eta = \frac{1}{2}[\xi, \eta].$$

Next, we define the Hopf fibering for a pseudosphere $S_1^2$. It is easy to see that the $\mathrm{Ad}(G)$-orbit of $\mathbf{k}' \in \mathfrak{g}$ is the pseudosphere:

$$S_1^2 = \{\xi \in \mathbb{E}_1^3 \mid \langle \xi, \xi \rangle = 1\}.$$

The Ad-action of $G$ on $S_1^2$ is transitive and isometric. The isotropy subgroup of $G$ at $\mathbf{k}'$ is the indefinite orthogonal group $\mathrm{SO}_1(2) = \{\xi_0 \mathbf{1} + \xi_3 \mathbf{k}' \mid \xi_0^2 - \xi_3^2 = 1\}$.

The Lie group $SO_1(2)$ is the hyperbola $H_0^1$ in a Minkowski plane $\mathbb{E}_1^2(\xi_0, \xi_3)$. (This is a Lorentz analogue of $S^1 \subset \mathbb{E}^2(\xi_0, \xi_1)$.) Note that the group $H_0^1$ is isomorphic to the multiplicative group $\mathbb{R}^\times$. (See Subsection 1.1.)

The natural projection $\pi \colon G \to S_1^2$, given by $\pi(g) = \mathrm{Ad}(g)\mathbf{k}'$ for all $g \in G$, defines a principal $H_0^1$-bundle over $S_1^2$. We call this fibering the *Hopf fibering of* $S_1^2$.

Denote the isotropy subgroup $H_0^1$ at $\mathbf{k}'$ by $K$ and its Lie algebra $\mathfrak{k} = \mathbb{R}\mathbf{k}'$. The tangent space of $S_1^2$ at the origin $\mathbf{k}'$ is given by $\mathfrak{m} = \mathbb{R}\mathbf{i} \oplus \mathbb{R}\mathbf{j}'$.

Let $\sigma$ be an involution of $\mathfrak{g}$ defined by $\sigma = \mathrm{Ad}(\mathbf{i}) = \Pi_\mathfrak{k} - \Pi_\mathfrak{m}$, where $\Pi_\mathfrak{k}$ and $\Pi_\mathfrak{m}$ are the projections from $\mathfrak{g}$ onto $\mathfrak{k}$ and $\mathfrak{m}$ respectively. The pair $(\mathfrak{g}, \sigma)$ is a symmetric Lie algebra data for the Lorentz symmetric space $S_1^2 = G/K$.

*Remark 1.7.* The Hopf fibering $\pi \colon H_1^3 \to S_1^2$ is naturally generalized to the principal $\mathbb{R}^\times$-bundle $\pi \colon H_n^{2n+1} \to \mathbb{C}'P_n$ over the paracomplex projective space $\mathbb{C}'P_n$ (cf. [22]).

## 2. Split-Complex Manifolds

Let $M$ be an $n$-manifold and $L(M)$ the linear frame bundle of $M$. A reduction of the structure group $\mathrm{GL}_n\mathbb{R}$ to

$$G = \left\{ \begin{pmatrix} A & O \\ O & B \end{pmatrix} \mid A \in \mathrm{GL}_{n_1}\mathbb{R},\ B \in \mathrm{GL}_{n_2}\mathbb{R} \right\}, \quad n_1 + n_2 = n$$

is called an *almost product structure* on $M$. Namely an almost product structure is a sub-bundle $P$ of $L(M)$ over $M$ with structure group $G$. (See Libermann [44].) This structure is an example of (a so-called) $G$-structure of infinite type. (See [39].)

An almost product structure corresponds to a tensor field $F$ of type $(1, 1)$ on $M$ such that $F^2 = I$.

LEMMA 2.1. *Let $M$ be an $n$-manifold. An almost product structure on $M$ corresponds to a pair of complementary distributions $(\mathcal{S}_1, \mathcal{S}_2)$ of rank $n_1$ and $n_2$ ($n = n_1 + n_2$).*

DEFINITION 2.2. A manifold $M$ together with an almost product structure is called an *almost product manifold*.

DEFINITION 2.3 (See [39]). An almost product structure $P$ is said to be *integrable* if every point of $M$ has a coordinate neighborhood $(U; x^1, \ldots, x^n)$ such that the cross section

$$\left( \frac{\partial}{\partial x^1}, \ldots, \frac{\partial}{\partial x^n} \right)$$

of $L(M)$ over $U$ is a cross section of $P$ over $U$.

LEMMA 2.4 (cf. [24]). *An almost product structure $(\mathcal{S}_1, \mathcal{S}_2)$ is integrable if and only if both distributions are integrable.*

Our interest is the case $n_1 = n_2 \geq 1$.

DEFINITION 2.5. An almost product structure $(\mathcal{S}_1, \mathcal{S}_2)$ with $n_1 = n_2 \geq 1$ is called an *almost split-complex structure*. In particular, an integrable almost split-complex structure is called a *split-complex structure*.

Let $M$ be a manifold with almost split-complex structure determined by the endomorphism field $K'$. Then the endomorphism field $K'$ induces a split-complex structure $K'_p$ on each tangent space $T_pM$ at $p \in M$. Thus every tangent space $T_pM$ is regarded as a free $\mathbb{C}'$-module. Let us denote the $\pm 1$-eigenspaces of $K'_p$ at $p$ by $T'_pM$ and $T''_pM$. Then the tangent bundle $TM$ splits into the following Whitney sum of real vector bundles:

$$TM = T_+M \oplus T_-M, \quad T_+M = \bigcup_{p \in M}(T_pM)_+, \quad T_-M = \bigcup_{p \in M}(T_pM)_-.$$

LEMMA 2.6. *Let $M$ be a 2-manifold and $K'$ an almost split-complex structure on $M$. Then $K'$ is integrable.*

*Proof.* Since $\dim M = 2$, $\operatorname{rank} T_+M = \operatorname{rank} T_-M = 1$. Hence both distributions

$$T_+M: p \mapsto (T_pM)_+, \qquad T_-M: p \mapsto (T_pM)_-$$

are integrable. □

*Remark 2.7.* A split-complex structure $K'$ is integrable if and only if its Nijenhuis torsion vanishes.

For general theory of $G$-structures, see Fujimoto [24] and Kobayashi [39].

The notion of holomorphic map between split-complex manifolds is introduced in the following way:

DEFINITION 2.8. Let $(M_1, K'_1)$, $(M_2, K'_2)$ be split-complex manifolds and $\phi: M_1 \to M_2$ a smooth map. Then $\phi$ is said to be a *holomorphic map* (or *paraholomorphic*) if $d\phi$ commutes with the split-complex structures:

$$d\phi \circ K'_1 = K'_2 \circ d\phi.$$

Similarly $\phi$ is said to be an *anti holomorphic map* (or *anti paraholomorphic*) if $d\phi$ satisfies

$$d\phi \circ K'_1 = -K'_2 \circ d\phi.$$

DEFINITION 2.9. Let $(M, K')$ be an almost split-complex manifold.

A semi-Riemannian metric $g$ is said to be a *split-Hermitian metric* if

$$g(K'X, K'Y) = -g(X, Y), \quad X, Y \in \mathfrak{X}(M).$$

An almost split-complex manifold together with split-Hermitian metric is called an almost *split-Hermitian manifold*.

On an almost split-Hermitian manifold $(M, K', g)$, the following tensor field $\Omega$:

$$\Omega(X, Y) := g(K'X, Y), \quad X, Y \in \mathfrak{X}(M)$$

is a 2-form on $M$. The 2-form $\Omega$ is called the *fundamental 2-form* of $(M, K', g)$.

DEFINITION 2.10. Let $(M, K', g)$ be an almost split-Hermitian manifold. Then $(M, K', g)$ is said to be a split-Kähler manifold (or a *para-Kähler manifold*) if $K'$ is integrable and $\Omega$ is closed.

Note that the fundamental 2-form of a split-Kähler manifold is symplectic.

On a 2-dimensional almost split-complex manifold, all the split-Hermitian metrics are split-Kähler. Hence, every 2-dimensional almost split-Hermitian manifold is split-Kähler. The following is proved in a way that is similar to the case of Kähler manifolds.

PROPOSITION 2.11. *Let $(M, K', g)$ be an almost split-Hermitian manifold. Then $(M, K', g)$ is a split-Kähler manifold if and only if the almost split-complex structure $K'$ is parallel with respect to the Levi-Civita connection of $g$.*

For general theory of split-Kähler manifolds, see [11]. In particular, for homogeneous split-Kähler manifolds, see Kaneyuki [35] and Kaneyuki and Kozai [36].

*Remark 2.12.* A split-Kähler manifold is said to be a *homogeneous split-Kähler manifold* if there exists a Lie group $G$ of holomorphic isometries which acts transitively on it.

As we will see later, the pseudosphere (Lorentz sphere) $S_1^2$ is a typical example of homogeneous split-Kähler manifold. In 2-dimensional geometry, there are three kinds of model geometries: Euclidean plane geometry, spherical geometry and hyperbolic plane geometry. In [56], H. Poincaré suggested the "4th-model geometry" in 2-dimension. The 4th-model geometry in his sense is the geometry of $S_1^2$. Kaneyuki gave a detailed exposition of the 4th-model geometry in [35].

## 3. Lorentz Surfaces

### 3.1. LORENTZ CONFORMAL STRUCTURE

In this subsection, we recall some basic facts on Lorentz conformal structures from the textbook [67].

Let $M$ be an oriented 2-manifold and $h_1$ and $h_2$ Lorentzian metrics on $M$. Then $h_1$ and $h_2$ are said to be *conformally equivalent* if there exists a smooth positive function $\mu$ on $M$ such that $h_2 = \mu h_1$. An equivalence class $\mathcal{C}$ of a Lorentzian metric on $M$ is called a *Lorentz conformal structure*. An orderd pair $(M, \mathcal{C})$ consisiting of an oriented surface and a Lorentz conformal structure compatible to the given orientation is called a *Lorentz surface*.

A local coordinate system $(u, v)$ is said to be *null* if

$$h\left(\frac{\partial}{\partial u}, \frac{\partial}{\partial u}\right) = h\left(\frac{\partial}{\partial v}, \frac{\partial}{\partial v}\right) = 0$$

for any $h \in \mathcal{C}$.

PROPOSITION 3.1 (cf. p. 13, [67]). *Let $(M, \mathcal{C})$ be a Lorentz surface. Then there exists, in some neighborhood of any point, a null coordinate system $(u, v)$.*

EXAMPLE 3.2. Let $\mathbb{E}_1^2$ be the Minkowski plane:

$$\mathbb{E}_1^2 = (\mathbb{R}^2, -\mathrm{d}x^2 + \mathrm{d}y^2)$$

with the natural Lorentz metric. The Lorentz surface determined by $\mathbb{E}_1^2$ is called the *conformal Minkowski plane*. It is easy to see that $u = x + y$, $v = -x + y$ give a global null coordinate system.

Next we recall the notion of a "box structure". Let $M$ be a 2-manifold with atlas $\mathcal{A}$. Two intersecting charts $(\mathcal{U}_1, \chi_1)$ and $(\mathcal{U}_2, \chi_2)$ in $\mathcal{A}$ are said to be $C^\square$-*related* if and only if the map $\chi_2 \circ \chi_1^{-1}$ taking $(u_1, v_1)$ to $(u_2, v_2)$ is given in some neighborhood of any point in $\chi_1(\mathcal{U}_1 \cap \mathcal{U}_2)$ by functions $u_2 = f(u_1)$ and $v_2 = g(v_1)$ with $f_{u_1} g_{v_1} > 0$. A $C^\square$-atlas $\mathcal{A}^\square$ on $M$ is a subatlas of $\mathcal{A}$ within which any two charts are $C^\square$-related.

A *box-structure* on $M$ is a maximal $C^\square$ atlas on $M$.

It is known that there is a natural one to one correspondence between box structures on a 2-manifold $M$ and Lorentz conformal structures on $M$.

THEOREM 3.3 ([67], p. 20). *There is a natural bijective correspondence between box structures $\mathcal{A}^\square$ and Lorentz conformal structures $\mathcal{C}$ on $M$.*

This fundamental fact motivates us to introduce the notion of "Lorentz holomorphic function" in Subsection 3.3.

3.2. SPLIT-KÄHLER STRUCTURE ON LORENTZ SURFACES

Let $(M, \mathcal{C})$ be a Lorentz surface. Take a null coordinate system $(u, v)$. Then the tensor field $K'$

$$K'\frac{\partial}{\partial u} = \frac{\partial}{\partial u}, \qquad K'\frac{\partial}{\partial v} = -\frac{\partial}{\partial v}$$

is defined on $M$ independent of the choice of null coordinates. Obviously $K'$ is an almost split-complex structure. Since $M$ is 2-dimensional, this almost split-complex structure is integrable. This structure $K'$ is called the *associated split-complex structure* of $(M, \mathcal{C})$. For any metric $h \in \mathcal{C}$, $(M, K', h)$ is split-Kähler.

LEMMA 3.4. *Let $(M, \mathcal{C})$ be a Lorentz surface. Then there exists a split-complex structure $K'$ on $M$.*

*Conversely, let $K'$ be an almost split-complex structure on $M$, then there exists a Lorentz conformal structure $\mathcal{C}$ on $M$ such that $K'$ is associated to $\mathcal{C}$.*

*Remark 3.5.* This lemma can be explained in terms of $G$-structures as follows: Let us denote by $\mathrm{CO}_1(2)$ the *linear conformal group* of $\mathbb{E}_1^2$:

$$\mathrm{CO}_1(2) := \big\{ A \in \mathrm{GL}_2\mathbb{R} \mid \langle A\mathbf{u}, A\mathbf{v} \rangle = c\langle \mathbf{u}, \mathbf{v} \rangle, \text{ for some } c > 0 \big\}.$$

Then one can see that $\mathrm{CO}_1(2)$ is isomorphic to the identity component of the following group:

$$\mathrm{GL}_1\mathbb{C}' := \big\{ A \in \mathrm{GL}_2\mathbb{R} \mid A\mathcal{K}' = \mathcal{K}'A \big\}.$$

Here $\mathcal{K}'$ is the natural split-complex structure of $\mathbb{E}_1^2$.

PROPOSITION 3.6. *Let $M_1$, $M_2$ be Lorentz surfaces and $\phi\colon M_1 \to M_2$ a smooth map. Take a null coordinate system $(u, v)$ on $M_1$ and $(f, g)$ on $M_2$. Then $\phi = (f(u, v), g(u, v))$ is holomorphic if and only if*

$$\frac{\partial f}{\partial v} = 0, \qquad \frac{\partial g}{\partial u} = 0.$$

*Analogously, $\phi = (f(u, v), g(u, v))$ is anti holomorphic if and only if*

$$\frac{\partial f}{\partial u} = 0, \qquad \frac{\partial g}{\partial v} = 0.$$

COROLLARY 3.7. *Let $(M_1, \mathcal{C}_1)$ and $(M_2, \mathcal{C}_2)$ be Lorentz surfaces and $\phi\colon M_1 \to M_2$ a smooth map. Take Lorentz isothermal coordinate systems $(x, y)$ and $(\xi, \eta)$ with respect to $\mathcal{C}_1$ and $\mathcal{C}_2$, respectively. Then $\phi = (\xi(x, y), \eta(x, y))$ is holomorphic if and only if*

$$\frac{\partial \xi}{\partial x} = \frac{\partial \eta}{\partial y}, \qquad \frac{\partial \xi}{\partial y} = \frac{\partial \eta}{\partial x}. \tag{3.1}$$

The system of partial differential equations in the corollary is regarded as *hyperbolic Cauchy–Riemann equations*.

*Remark 3.8.* Let $\phi\colon D \subset \mathbb{C}' \to \mathbb{C}'$ be a map on a region of $\mathbb{C}'$. Then $\phi$ is said to be *split-complex differentiable* (or *paracomplex differentiable*) at $\zeta \in D$ if there exists a split-complex number $\alpha(\zeta)$ such that

$$\phi(\zeta + \mathrm{h}) - \phi(\zeta) = \alpha(\zeta)\mathrm{h} + \sigma(\zeta, \mathrm{h})(h_1^2 + h_2^2),$$

where $h = h_1 + h_2 k'$ and $\sigma(\zeta, h) \to 0$ as $h \to 0$. It is known that a smooth map $\phi \colon D \subset \mathbb{C}' \to \mathbb{C}'$ is split-complex differentiable in $D$ if and only if $\phi$ is Lorentz holomorphic. (See [11, 21].)

In [11, 27], the notion of "paracomplex analytic function" is introduced. A map $\phi \colon \mathbb{C}' \to \mathbb{C}'$ is said to be *paracomplex analytic* if it has a power series expansion in $\mathbb{C}'$. The paracomplex analyticity implies Lorentz holomorphicity. But the converse is not always true.

*Remark 3.9.* Let $(M, \mathcal{C})$ be a Lorentz surface. Take a Lorentzian metric $g \in \mathcal{C}$. Denote by $\Box_g$ the Laplace–Bertrami operator of $g$ (frequently called *d'Alembert operator*). With respect to an Lorentz isothermal coordinate system $(x, y)$ such that $g = e^\omega(-dx^2 + dy^2)$, $\Box_g$ is given by

$$\Box_g = e^{-\omega}\left(-\frac{\partial^2}{\partial x^2} + \frac{\partial^2}{\partial y^2}\right).$$

A smooth function $f \in C^\infty(M)$ is said to be a (*Lorentz*) *harmonic function* if $\Box f = 0$. The explicit formula of $\Box$ implies that the *harmonicity of functions* is invariant under conformal transformations of $(M, \mathcal{C})$. The hyperbolic Cauchy–Riemann equations imply that each component of a holomorphic map between Lorentz surfaces is a (Lorentz) harmonic function. (See Definition 3.13.)

### 3.3. REAL HOLOMORPHIC FUNCTIONS

Let $M$ be a Lorentz surface. Denote by $K'$ the associated split-complex structure of $M$. Then $TM$ splits into

$$TM = T_+M \oplus T_-M.$$

According to this splitting, $T^*M$ splits into

$$T^*M = T'M \oplus T''M.$$

Denote by $A^{1,0}(M)$ and $A^{0,1}(M)$ the linear space of all smooth sections of $T'M$ and $T''M$ respectively.

By definition, every element $A \in A^{1,0}(M)$ has the following local expression

$$A = a\, du.$$

Similarly, every element $B \in A^{0,1}(M)$ has the following local expression

$$B = b\, dv.$$

Let us denote the space of all smooth tensor fields of type $(0, 2)$ by $\mathcal{T}^0_2(M)$. We introduce the following subspaces of $\mathcal{T}^0_2(M)$:

$$\mathcal{T}^{2,0}(M) = \Gamma(T'M \otimes T'M),$$

$$\mathcal{T}^{1,1}(M) = \Gamma(T'M \otimes T''M),$$

$$\mathcal{T}^{0,2}(M) = \Gamma(T''M \otimes T''M).$$

Denote by $A^2(M)$ and $S^2(M)$ the space of all smooth differential 2-forms and space of all smooth symmetric 2-tensor fields on $M$ respectively:

$$S^2(M) = \Gamma(T^*M \odot T^*M), \qquad A^2(M) = \Gamma(T^*M \wedge T^*M).$$

Put $S^{p,q}(M) := S^2(M) \cap \mathcal{T}^{p,q}(M)$ and $A^{p,q}(M) := A^2(M) \cap \mathcal{T}^{p,q}(M)$. Then $S^2(M)$ splits into

$$S^2(M) = S^{2,0}(M) \oplus S^{1,1}(M) \oplus S^{0,2}(M),$$

$$A^2(M) = A^{1,1}(M), \qquad A^{2,0}(M) = A^{0,2}(M) = \{0\}.$$

The differential operator d is decomposed with respect to the conformal structure:

$$\mathrm{d} = \mathrm{d}' + \mathrm{d}'', \quad \mathrm{d}' := \frac{\partial}{\partial u}\mathrm{d}u, \ \mathrm{d}'' := \frac{\partial}{\partial v}\mathrm{d}v.$$

It is easy to see that $\mathrm{d}'$ and $\mathrm{d}''$ are independent of the choice of null coordinate system. These operators are regarded as exterior differential operators:

$$\mathrm{d}' \colon A^0(M) \to A^{1,0}(M), \qquad \mathrm{d}'' \colon A^0(M) \to A^{0,1}(M), \quad A^0(M) := C^\infty(M).$$

These operators can be naturally extended to

$$\mathrm{d}' \colon A^{0,1}(M) \to A^{1,1}(M), \qquad \mathrm{d}'' \colon A^{1,0}(M) \to A^{1,1}(M).$$

DEFINITION 3.10. Let $f \colon M \to \mathbb{R}$ be a smooth function. Then $f$ is said to be a (real) *Lorentz holomorphic function* if $\mathrm{d}''f = 0$. Similarly, $f$ is said to be a (real) *Lorentz anti holomorphic function* if $\mathrm{d}'f = 0$.

It is easy to see that $f \in C^2(M)$ is Lorentz harmonic if and only if $\mathrm{d}''\mathrm{d}'f = 0$. In particular, (Lorentz) $\pm$-holomorphic functions are harmonic.

PROPOSITION 3.11 (d'Alembert formula). *Let $M$ be a Lorentz surface and $h$ a Lorentz harmonic function. Then there exist a Lorentz holomorphic function $f$ and a Lorentz anti holomorphic function $g$ such that $h = f + g$.*

Next we introduce the notion of *holomorphicity* of 1-forms.

DEFINITION 3.12. Let $A$ be a $(1, 0)$-type 1-form on $M$. Then $A$ is said to be a *Lorentz holomorphic 1-form* if $\mathrm{d}''A = 0$. Similarly, a $(0, 1)$-type 1-form $B$ is said to be a *Lorentz anti holomorphic 1-form* if $\mathrm{d}'B = 0$.

The induced split-complex structure on $T^*M$ may be considered as a Lorentzian analogue of Hodge star operator. Hereafter, we use the notation $*$ for $K'\colon T^*M \to T^*M$ and call it *star operator*. The star operator has the following local expression:

$$*\mathrm{d}u = \mathrm{d}u, \qquad *\mathrm{d}v = -\mathrm{d}v.$$

Take a Lorentzian metric $g \in \mathcal{C}$. Express the metric $g$ as $g = e^{\omega}\,\mathrm{d}u\,\mathrm{d}v$. The area element $\mathrm{d}A_g$ of $(M, g)$ is $\mathrm{d}A_g = e^{\omega}\mathrm{d}u \wedge \mathrm{d}v/2$. The *co-differential operator* $\delta = \delta_g$ with respect to $g$ is defined by

$$\delta := * \circ \mathrm{d} \circ *\colon A^p(M) \to A^{p-1}(M), \quad *1 = \mathrm{d}A_g,\ *\mathrm{d}A_g = 1.$$

The Laplace–Bertrami operator (d'Alembert operator) $\Box = \Box_g\colon A^0(M) \to A^0(M)$ of $(M, g)$ is expressed by $\Box_g = \delta_g\,\mathrm{d}$.

DEFINITION 3.13. Let $(M, g)$ be a Lorentzian 2-manifold. The *d'Alembert operator* $\Box_g$ of $(M, g)$ acting on $A^p(M)$ is defined by

$$\Box_g := \mathrm{d}\delta_g + \delta_g\,\mathrm{d}\colon A^p(M) \to A^p(M), \quad p = 0, 1, 2.$$

A $p$-form $\theta$ is said to be a *Lorentz harmonic $p$-form* if $\Box_g\theta = 0$.

COROLLARY 3.14. *Let $\theta$ be a 1-form on a Lorentz surface $M$. If $\theta$ satisfies*

$$\mathrm{d}\theta = 0, \qquad \mathrm{d} * \theta = 0,$$

*then $\theta$ is Lorentz harmonic.*

Note that harmonicity depends only on the conformal structure. Namely harmonicity for 1-forms makes sense over Lorentz surfaces.

PROPOSITION 3.15. *Let $f\colon M \to \mathbb{R}$ a smooth function. Then $f$ is Lorentz harmonic if and only if*

$$\mathrm{d} * \mathrm{d}f = 0.$$

*Remark 3.16.* Some authors call the d'Alembert operator *wave operator*. Harmonic maps of Lorentz surfaces into semi-Riemannian manifolds are sometimes called *wave maps*.

DEFINITION 3.17. Let $(M_1, \mathcal{C}_1)$ and $(M_2, \mathcal{C}_2)$ be Lorentz surfaces. A smooth map $\phi\colon M_1 \to M_2$ is said to be a *Lorentz-conformal map* if it preserves the conformal structure, i.e.,

$$\phi^* h_2 \in \mathcal{C}_1$$

for any $h_2 \in \mathcal{C}_2$.

Proposition 3.6 and Corollary 3.7 imply that Lorentz-conformality of $\phi\colon M_1 \to M_2$ is equivalent to the (para) holomorphicity of $\phi$.

EXAMPLE 3.18 (Stereographic projection). Let $S_1^2$ be the unit pseudosphere in $\mathbb{E}_1^3$:

$$S_1^2 = \{(x_1, x_2, x_3) \in \mathbb{E}_1^3 \mid -x_1^2 + x_2^2 + x_3^2 = 1\}.$$

We identify the Minkowski plane $\mathbb{E}_1^2$ with the timelike plane $x_3 = 0$ in $\mathbb{E}_1^3$. Let $\mathcal{N} = (0, 0, -1)$ and $\mathcal{S} = (0, 0, 1)$ be the *north pole* and *south pole* of $S_1^2$. The *stereographic projection* $\mathcal{P}_+\colon S_1^2 \setminus \{x_3 = -1\} \to \mathbb{E}_1^2 \setminus H_0^1$ with respect to the north pole is defined as a map which assigns a point $p$ to the intersection of $\mathbb{E}_1^2$ and the line through $P$ and $\mathcal{N}$. The stereographic projection $\mathcal{P}_+$ is given explicitly by

$$\mathcal{P}_+(x_1, x_2, x_3) = \left(\frac{x_1}{x_3 + 1}, \frac{x_2}{x_3 + 1}, 0\right).$$

This explicit representation implies that $\mathcal{P}_+$ is a Lorentz-conformal map from $S_1^2 \setminus \{x_3 = -1\}$ onto $\mathbb{E}_1^2 \setminus \{-x_1^2 + x_2^2 = -1\}$. The inverse map $\mathcal{P}_+^{-1}$ is

$$\mathcal{P}_+^{-1}(u, v) = \frac{1}{-u^2 + v^2 + 1}(2u, 2v, u^2 - v^2 + 1).$$

The stereographic projection $\mathcal{P}_-$ with respect to the south pole $\mathcal{S}$ is analogously defined.

## Part II. Timelike Immersions

## 4. Timelike Surfaces

Let $M$ be a connected 2-manifold and $\varphi\colon M \to \mathbb{E}_1^3$ an immersion. The immersion $\varphi$ is said to be *timelike* if the induced metric I of $M$ is Lorentzian. Hereafter we may assume that $M$ is an orientable timelike surface in $\mathbb{E}_1^3$ (immersed by $\varphi$). It is worthwhile to remark that there exists no compact timelike surface in $\mathbb{E}_1^3$. (See [55], p. 125.)

The induced Lorentzian metric I determines a Lorentz conformal structure $\mathcal{C}_I$ on $M$. With respect to $\mathcal{C}_I$, we regard $M$ as a Lorentz surfaces and $\varphi$ a conformal immersion.

On a timelike surface $M$, there exists a local null coordinate system $(u, v)$ with respect to the conformal structure $\mathcal{C}_I$ (Proposition 3.1),

$$\mathrm{I} = e^\omega \, du \, dv.$$

Define a local coordinate system $(x, y)$ by $x = (u - v)/2$, $y = (u + v)/2$. Then the first fundamental form I is written as

$$\mathrm{I} = e^\omega(-dx^2 + dy^2). \tag{4.1}$$

The local coordinate system $(x, y)$ is called the *Lorentz isothermal coordinate system* associated to $(u, v)$.

Now, let $N$ be a unit normal vector field to $M$. The vector field $N$ is spacelike, i.e., $\langle N, N \rangle = 1$, since $M$ is timelike. The vector fields $\varphi_u, \varphi_v$ as well as the normal $N$ define a moving frame. We shall define a moving frame $\mathfrak{s}$ by $\mathfrak{s} = (\varphi_u, \varphi_v, N)$. The moving frame $\mathfrak{s}$ satisfies the following Frenet (or Gauss–Weingarten) equations:

$$\mathfrak{s}_u = \mathfrak{s}\mathcal{U}, \qquad \mathfrak{s}_v = \mathfrak{s}\mathcal{V}, \tag{4.2}$$

$$\mathcal{U} = \begin{pmatrix} \omega_u & 0 & -H \\ 0 & 0 & -2e^{-\omega}Q \\ Q & \frac{1}{2}He^{\omega} & 0 \end{pmatrix}, \quad \mathcal{V} = \begin{pmatrix} 0 & 0 & -2e^{-\omega}R \\ 0 & \omega_v & -H \\ \frac{1}{2}He^{\omega} & R & 0 \end{pmatrix}, \tag{4.3}$$

where $Q := \langle \varphi_{uu}, N \rangle$, $R := \langle \varphi_{vv}, N \rangle$, $H = 2e^{-\omega}\langle \varphi_{uv}, N \rangle$. It is easy to see that $Q^{\#} := Q\, du^2$ and $R^{\#} := R\, dv^2$ are globally defined on $(M, \mathcal{C}_l)$. Moreover $Q^{\#} \in S^{2,0}(M)$ and $R^{\#} \in S^{0,2}(M)$. We call the differential $Q^{\#} + R^{\#}$, the *Hopf differential* of $M$.

It is easy to see that $H$ is the mean curvature of $M$. The *second fundamental form* $\mathbb{I}$ of $M$ defined by

$$\mathbb{I} = -\langle d\varphi, dN \rangle \tag{4.4}$$

is described relative to the Lorentz isothermal coordinates $(x, y)$ as follows:

$$\mathbb{I} = \begin{pmatrix} Q + R - He^{\omega} & Q - R \\ Q - R & Q + R + He^{\omega} \end{pmatrix}. \tag{4.5}$$

The Gaussian curvature $K$ of $M$ is given by

$$K = \det(\mathbb{I} \cdot \mathrm{I}^{-1}).$$

(See [55], p. 107.) The second fundamental form $\mathbb{I}$ is decomposed as

$$\mathbb{I} = Q^{\#} + H\mathrm{I} + R^{\#}$$

along the decomposition

$$S^2(M) = S^{2,0}(M) \oplus S^{1,1}(M) \oplus S^{0,2}(M).$$

The *Gauss equation* which describes a relation between $K$, $H$ and $Q$ takes the following form:

$$H^2 - K = 4e^{-2\omega}QR.$$

It is easy to see that the zeros of Hopf differential coincide with the umbilic points of $M$.

The Gauss–Codazzi equation, i.e. the integrability condition of the Frenet equation,

$$\mathcal{V}_u - \mathcal{U}_v + [\mathcal{U}, \mathcal{V}] = 0$$

has the following form:

$$\omega_{uv} + \frac{1}{2}H^2 e^\omega - 2QR e^{-\omega} = 0, \tag{4.6}$$

$$H_u = 2e^{-\omega} Q_v, \tag{4.7}$$

$$H_v = 2e^{-\omega} R_u. \tag{4.8}$$

Equations (4.7)–(4.8) show that the constancy of the mean curvature $H$ is equivalent to the condition $Q_v = R_u = 0$. As in the case of spacelike surfaces, we have the following Bonnet-type theorem.

PROPOSITION 4.1. *Every timelike constant mean curvature surface which is not totally umbilic has a local one parameter family of isometric deformations preserving the mean curvature.*

*Proof.* Take a simply connected null coordinate region $(\mathbb{D}; u, v)$ of a timelike constant mean curvature surface. Then the Gauss–Codazzi equations (4.7)–(4.8) are invariant under the deformation:

$$Q \longmapsto Q_\lambda := \lambda Q, \qquad R \longmapsto R_\lambda := \lambda^{-1} R, \quad \lambda \in \mathbb{R}^\times. \tag{4.9}$$

Integrating the Gauss–Weingarten equations (4.2) with deformed data $Q_\lambda, R_\lambda$ over $\mathbb{D}$, one obtains a one-parameter family of timelike surfaces $\{\varphi_\lambda\}$ defined on $\mathbb{D}$. The deformation (4.9) does not effect the induced metric and the mean curvature. Hence all the surfaces $\{\varphi_\lambda\}$ are isometric and have the same constant mean curvature. □

The one-parameter family $\{\varphi_\lambda\}$ is called the *associated family* of $\varphi$.

*Remark 4.2.* A non flat totally umbilic timelike surface is congruent to an open portion of a pseudosphere $S_1^2(r) = \{\xi \in \mathbb{E}_1^3 \mid \langle \xi, \xi \rangle = r^2\}$ of radius $r > 0$. (Hence if $M$ is complete, $M$ is congruent to $S_1^2(r)$.) Note that there is no totally umbilic timelike surface of (constant) negative curvature. (See [55], p. 116.)

*Remark 4.3.* Since the metric I is indefinite, the condition $QR = 0$ is not equivalent to $Q = R = 0$. This fact yields that even if the eigenvalues of the *shape operator* $S := -dN$ repeat, then $M$ is not necessarily totally umbilic [51].

### 4.1. GAUSS MAPS

Next, we shall define the Gauss map of a timelike surface. Let $M$ be a timelike surface and $N$ a unit normal vector field to $M$. For each $p \in M$ the point $\psi(p)$ of $\mathbb{E}_1^3$ canonically corresponding to the vector $N_p$ lies in a unit pseudosphere since $N$ is spacelike. The resulting smooth mapping $\psi\colon M \to S_1^2$ is called the *Gauss map* of $M$. The constancy of mean or Gaussian curvature is characterized by the harmonicity of the Gauss map (see [67]).

PROPOSITION 4.4. *The Gauss map of a timelike surface is harmonic if and only if the mean curvature is constant.*

PROPOSITION 4.5. *Let $M$ be a timelike surface. Assume that the Gaussian curvature $K$ is nowhere zero on $M$ and has a constant sign. Then the second fundamental form $\mathbb{II}$ gives $M$ another (semi-) Riemannian metric. With respect to this metric $\mathbb{II}$, the Gauss map of $M$ is harmonic if and only if $K$ is constant.*

The *third fundamental form* $\mathbb{III}$ of $(M, \varphi)$ is defined by

$$\mathbb{III}(X, Y) := \langle SX, SY \rangle.$$

Here $S$ is the shape operator of $M$. The third fundamental form satisfies

$$\mathbb{III} = 2H\mathbb{II} - K\mathrm{I}.$$

Direct computation shows

$$\langle \psi_* X, \psi_* Y \rangle = -K\mathrm{I}(X, Y).$$

Thus we get the following.

PROPOSITION 4.6. *Let $\varphi\colon M \to \mathbb{E}_1^3$ be a timelike surface with Gauss map $\psi\colon M \to S_1^2$. Away from flat points, the conformality of $\psi$ is equivalent to the condition $H = 0$.*

### 4.2. LAX FORMS

Now, we recall the $2 \times 2$ matrix-representation of the Gauss–Codazzi equations obtained by the first named author [34]. Let $\varphi\colon M \to \mathbb{E}_1^3$ be a timelike surface as in the preceding section. Since $M$ is orientable, there exists a global unit normal vector field $N$ to $M$ up to sign. Hereafter we choose the following sign convention. Take a local Lorentz isothermal coordinate system $(x, y)$ which is compatible to the orientation. Then $N$ is locally expressed as

$$N = (\varphi_x \times \varphi_y) / |\varphi_x \times \varphi_y|.$$

Define a matrix-valued function $\Phi$ by

$$\mathrm{Ad}(\Phi)(\mathbf{i}, \mathbf{j}', \mathbf{k}') = (\mathrm{e}^{-\frac{\omega}{2}}\varphi_x,\ \mathrm{e}^{-\frac{\omega}{2}}\varphi_y,\ N). \tag{4.10}$$

To derive the linear differential equations corresponding to (4.2), we introduce the following $\mathbb{H}'$-valued functions $U$ and $V$.

$$U = \Phi^{-1}\Phi_u, \qquad V = \Phi^{-1}\Phi_v.$$

Using the compatibility condition $\varphi_{uv} = \varphi_{vu}$ for (4.10) and the Gauss–Weingarten formulae (4.2), we get the following (Lax) pair:

$$U = u_0 \begin{pmatrix} 1 & 0 \\ 0 & 1 \end{pmatrix} + \begin{pmatrix} -\frac{1}{4}\omega_u & -Q e^{-\frac{\omega}{2}} \\ \frac{H}{2}e^{\frac{\omega}{2}} & \frac{1}{4}\omega_u \end{pmatrix}, \tag{4.11}$$

$$V = v_0 \begin{pmatrix} 1 & 0 \\ 0 & 1 \end{pmatrix} + \begin{pmatrix} \frac{1}{4}\omega_v & -\frac{H}{2}e^{\frac{\omega}{2}} \\ R e^{-\frac{\omega}{2}} & -\frac{1}{4}\omega_v \end{pmatrix} \tag{4.12}$$

under the identification (1.3). Only the coefficients $u_0$ and $v_0$ of $\mathbf{1}$ are still not determined. Recall that $\Phi$ was defined by (4.10) up to a multiplication by a scalar factor ($SO_1(2)$-valued function). Direct calculations show that the compatibility condition of the above linear differential equation is equivalent to the Gauss–Codazzi equations (4.6)–(4.8) if and only if $\frac{\partial}{\partial u} v_0 = \frac{\partial}{\partial v} u_0$ holds.

PROPOSITION 4.7 ([34]). *Under the identification (1.3), the moving frame* $\mathfrak{s} = (e^{-\frac{\omega}{2}}\varphi_x, e^{-\frac{\omega}{2}}\varphi_y, N)$ *of a timelike surface is described by the formula (4.10). Here* $\Phi$ *is an* $\mathbb{H}'$*-valued function which satisfies*

$$\frac{\partial \Phi}{\partial u} = \Phi \begin{pmatrix} u_0 - \frac{1}{4}\omega_u & -Q e^{-\frac{\omega}{2}} \\ \frac{H}{2}e^{\frac{\omega}{2}} & u_0 + \frac{1}{4}\omega_u \end{pmatrix}, \qquad \frac{\partial \Phi}{\partial v} = \Phi \begin{pmatrix} v_0 + \frac{1}{4}\omega_v & -\frac{H}{2}e^{\frac{\omega}{2}} \\ R e^{-\frac{\omega}{2}} & v_0 - \frac{1}{4}\omega_v \end{pmatrix}.$$

*The entries* $u_0$ *and* $v_0$ *are defined by (4.10) up to a multiplication of a scalar factor such that* $\frac{\partial}{\partial u}v_0 = \frac{\partial}{\partial v}u_0$.

## 5. Weierstrass-Type Representation

Here we briefly review the Weierstrass representation for a harmonic map $\psi \colon \mathbb{D} \to S_1^2$ from a simply connected Lorentz surface into the unit pseudosphere $S_1^2$. For our approach, we realize $S_1^2$ as a hyperboloid of one sheet in Minkowski 3-space. The harmonic map equation for $\psi$ is

$$\Box \psi + \langle d\psi, d\psi \rangle \psi = 0. \tag{5.1}$$

Next, as before, we express $S_1^2 = G/K = SL_2\mathbb{R}/SO_1(2)$ as a Lorentz symmetric space. Take a lift $\Psi \colon \mathbb{D} \to G$. Any such lift is called a *framing* of $\psi$. Denote by $\alpha$ the pull-back 1-from of the Maurer–Cartan from of $G$:

$$\alpha = \Psi^{-1}d\Psi.$$

We decompose $\alpha$ as

$$\alpha = \alpha_{\mathfrak{k}} + \alpha_{\mathfrak{m}}$$

along the Lie algebra decomposition $\mathfrak{g} = \mathfrak{k} \oplus \mathfrak{m}$. The Lie subspace $\mathfrak{m}$ is the tangent space of $S_1^2$ at the origin $\mathbf{k}'$.

With respect to the conformal structure of $\mathbb{D}$, $\alpha$ is decomposed as

$$\alpha = \alpha_\mathfrak{k} + \alpha'_\mathfrak{m} + \alpha''_\mathfrak{m}.$$

Here $\alpha'_\mathfrak{m}$ and $\alpha''_\mathfrak{m}$ are the $(1, 0)$ and $(0, 1)$ parts of $\alpha_\mathfrak{m}$, respectively.

Now we insert a spectral parameter $\lambda \in \mathbb{R}^+$ by the rule:

$$\alpha_\lambda := \alpha_\mathfrak{k} + \lambda \alpha'_\mathfrak{m} + \lambda^{-1} \alpha''_\mathfrak{m}. \tag{5.2}$$

Then the harmonicity of $\psi$ is equivalent to the flatness of the connection $A_\lambda := \mathrm{d} + \alpha_\lambda$;

$$F(A_\lambda) = \mathrm{d}\alpha_\lambda + \frac{1}{2}[\alpha_\lambda \wedge \alpha_\lambda] = 0, \quad \lambda \in \mathbb{R}^+.$$

Thus we arrive at the following recipe constructing harmonic maps of $\mathbb{D}$ into $G/K$.

Find an $\mathbb{R}^+$-family of flat connections $A_\lambda$ whose $\lambda$-dependence verifies (5.2). Then by integrating $\mathrm{d}\Psi_\lambda = \Psi_\lambda \alpha_\lambda$, we obtain a family of framings $\Psi_\lambda \colon \mathbb{D} \times \mathbb{R}^+ \to G$. The one-parameter family $\Psi_\lambda$ is called an *extended framing*. By projecting these maps to $G/K$, we finally obtain a one-parameter family $\{\psi_\lambda\}$ of harmonic maps. The one-parameter family $\{\psi_\lambda\}$ is called the *associated family* of a harmonic map $\psi = \psi_1$.

For the construction of extended framings, twisted loop groups are useful.

Let us denote by $\Lambda G$ the *free loop group* of $G$:

$$\Lambda G = \{\gamma \colon S^1 \to G\}.$$

The twisted loop group $\Lambda G_\sigma$ is

$$\Lambda G_\sigma = \{\gamma \in \Lambda G \mid \sigma(\gamma(\lambda)) = \gamma(-\lambda)\}.$$

The corresponding Lie algebras are given by

$$\Lambda\mathfrak{g} = \{\xi \colon S^1 \to \mathfrak{g}\}, \qquad \Lambda\mathfrak{g}_\sigma = \{\xi \in \Lambda\mathfrak{g} \mid \sigma(\xi(\lambda)) = \xi(-\lambda)\}.$$

There are several possibilities of introducing topologies on these groups. To simplify the discussion, we present the topological arguments of loop groups and loop algebras in Appendix A.

Next, we introduce the following subgroups:

$$\Lambda^+ G_\sigma = \left\{\gamma = \sum_{j \geqslant 0} \gamma_j \lambda^j \in \Lambda G_\sigma\right\}, \qquad \Lambda^- G_\sigma = \left\{\gamma = \sum_{j \leqslant 0} \gamma_j \lambda^j \in \Lambda G_\sigma\right\},$$

$$\Lambda_*^+ G_\sigma = \left\{ \gamma = \mathbf{1} + \sum_{j>0} \gamma_j \lambda^j \in \Lambda G_\sigma \right\},$$

$$\Lambda_*^- G_\sigma = \left\{ \gamma = \mathbf{1} + \sum_{j<0} \gamma_j \lambda^j \in \Lambda G_\sigma \right\}.$$

The corresponding Lie algebras are denoted by $\Lambda^+ \mathfrak{g}_\sigma$, $\Lambda^- \mathfrak{g}_\sigma$, $\Lambda_*^+ \mathfrak{g}_\sigma$, $\Lambda_*^+ \mathfrak{g}_\sigma$ respectively.

Finally, let $\tilde{\Lambda} G$ be the subset of $\Lambda G$ whose elements, as maps defined on $S^1$, admit analytic continuations to $\mathbb{C}^\times$. Then $\tilde{\Lambda} G$ is a Lie subgroup of $\Lambda G$. Similarly, we define $\tilde{\Lambda} G_\sigma$, $\tilde{\Lambda}^\pm G_\sigma$ and $\tilde{\Lambda}_*^\pm G_\sigma$.

Let $\Psi_\lambda \colon \mathbb{D} \times \mathbb{R}^+ \to G$ be an extended framing. Then, since $\Psi_\lambda$ is analytic in $\lambda$, $\Psi_\lambda$ has an analytic continuation to $\mathbb{D} \times \mathbb{C}^\times$. We denote this continuation by the same letter. Then $\Psi_\lambda$ can be viewed as a map $\mathbb{D} \to \tilde{\Lambda} G_\sigma$.

Very important ingredients of the Weierstrass-type representation of harmonic maps (or timelike constant mean curvature surfaces) are the following decomposition theorems obtained in [17]:

THEOREM 5.1 (Birkhoff decomposition of $\tilde{\Lambda} G_\sigma$). *We have*

$$\tilde{\Lambda} G_\sigma = \bigsqcup_{w \in \mathcal{T}} \tilde{\Lambda}^- G_\sigma \cdot w \cdot \tilde{\Lambda}^+ G_\sigma.$$

*Here $\mathcal{T}$ denotes the group of homomorphisms from $S^1$ into the subgroup of diagonal matrices of $\mathrm{SL}_2 \mathbb{R}$. Moreover, the multiplication maps*

$$\tilde{\Lambda}_*^- G_\sigma \times \tilde{\Lambda}^+ G_\sigma \to \tilde{\Lambda} G_\sigma, \qquad \tilde{\Lambda}_*^+ G_\sigma \times \tilde{\Lambda}^- G_\sigma \to \tilde{\Lambda} G_\sigma$$

*are diffeomorphisms onto the open dense subsets $\mathcal{B}(-, +)$ and $\mathcal{B}(+, -)$ of $\tilde{\Lambda} G_\sigma$, respectively, called the big cells of $\tilde{\Lambda} G_\sigma$. In particular if $\gamma$ is an element of $\mathcal{B} = \mathcal{B}(-, +) \cap \mathcal{B}(+, -)$, then $\gamma$ has unique decompositions*:

$$\gamma = \gamma_- \cdot \ell_+ = \gamma_+ \cdot \ell_-, \quad \gamma_\pm \in \tilde{\Lambda}_*^\pm G_\sigma, \ \ell_\pm \in \tilde{\Lambda}^\pm G_\sigma,$$

*respectively.*

THEOREM 5.2 (Iwasawa decomposition of $\tilde{\Lambda} G_\sigma \times \tilde{\Lambda} G_\sigma$). *Let $\Delta(\tilde{\Lambda} G_\sigma \times \tilde{\Lambda} G_\sigma)$ denote the diagonal subgroup of $\tilde{\Lambda} G_\sigma \times \tilde{\Lambda} G_\sigma$;*

$$\Delta(\tilde{\Lambda} G_\sigma \times \tilde{\Lambda} G_\sigma) = \{(g, g) \in \tilde{\Lambda} G_\sigma \times \tilde{\Lambda} G_\sigma\}.$$

*Then we have*

$$\tilde{\Lambda} G_\sigma \times \tilde{\Lambda} G_\sigma = \bigsqcup \Delta(\tilde{\Lambda} G_\sigma \times \tilde{\Lambda} G_\sigma) \cdot (\mathbf{1}, w) \cdot (\tilde{\Lambda}^- G_\sigma \times \tilde{\Lambda}^+ G_\sigma),$$

*where $w \in \mathcal{T}$ is as in the Birkhoff Decomposition Theorem.*

*Moreover, the multiplication maps*

$$\Delta(\tilde{\Lambda}G_\sigma \times \tilde{\Lambda}G_\sigma) \times (\tilde{\Lambda}_*^- G_\sigma \times \tilde{\Lambda}^+ G_\sigma) \to \tilde{\Lambda}G_\sigma \times \tilde{\Lambda}G_\sigma,$$

$$\Delta(\tilde{\Lambda}G_\sigma \times \tilde{\Lambda}G_\sigma) \times (\Lambda_*^+ G_\sigma \times \tilde{\Lambda}^- G_\sigma) \to \tilde{\Lambda}G_\sigma \times \tilde{\Lambda}G_\sigma$$

*are diffeomorphisms onto the open dense subsets* $\mathcal{I}(+,-)$ *and* $\mathcal{I}(-,+)$ *of* $\tilde{\Lambda}G_\sigma \times \tilde{\Lambda}G_\sigma$ – *called the big cells of* $\tilde{\Lambda}G_\sigma \times \tilde{\Lambda}G_\sigma$.

These decomposition theorems can be deduced from the Birkhoff Decomposition Theorem for $\tilde{\Lambda}G$ given in Appendix. (For more details, see [17].)

Let $\Psi\colon \mathbb{D} \to \tilde{\Lambda}G_\sigma$ be an extended framing. Then, away from a discrete set $\mathcal{S} \subset \mathbb{D}$ (cf. [15, 16, 19]), $\Psi$ can be decomposed as

$$\Psi = \Psi_- L_+ = \Psi_+ L_-, \quad \Psi_\pm \in \tilde{\Lambda}_*^\pm G_\sigma, \; L_\pm \in \tilde{\Lambda}^\pm G_\sigma$$

according to the Birkhoff Decomposition Theorem. Then one can easily check that

$$\xi' := \Psi_+^{-1} d\Psi_+ = \lambda \begin{pmatrix} 0 & p_2'(u) \\ p_1'(u) & 0 \end{pmatrix} du,$$

$$\xi'' := \Psi_-^{-1} d\Psi_- = \lambda^{-1} \begin{pmatrix} 0 & -p_1''(v) \\ p_2''(v) & 0 \end{pmatrix} dv.$$

The pair of 1-forms $\{\xi', \xi''\}$ is called the *normalized potential* of the harmonic map $\psi$.

Conversely, let $\{\xi', \xi''\}$ be a normalized potential defined on a simply connected Lorentz surface $\mathbb{D}$ with global null coordinate system $(u, v)$. Without loss of generality, we may assume that $\mathbb{D}$ is a simply connected region of $\mathbb{E}_1^2$ which contains the origin $(0,0)$, since we assumed the existence of global null coordinate system (see [67]). Then one can solve the initial value problems:

$$\Psi_+^{-1} d'\Psi_+ = \xi', \qquad \Psi_-^{-1} d''\Psi_- = \xi'', \qquad \Psi_+(u=0) = \Psi_-(v=0) = \mathbf{1}.$$

Now we apply the Iwasawa decomposition Theorem to $(\Psi_+, \Psi_-)$:

$$(\Psi_+, \Psi_-) = (\Psi, \Psi)(L_-^{-1}, L_+^{-1}). \tag{5.3}$$

Then we have

THEOREM 5.3 ([17]). *The mapping $\Psi$ obtained by the Iwasawa splitting* (5.3) *is an extended framing defined on a neighborhood of the origin.*

Hence we get a one-parameter family of harmonic maps $\psi_\lambda = \mathrm{Ad}(\Psi)\mathbf{k}'$.

The construction of harmonic maps (or corresponding timelike constant mean curvature surfaces) from (normalized) potentials is called the *Weierstrass representation* of harmonic maps, *DPW-construction* or *nonlinear d'Alembert formula* [40].

This recipe can be generalized to harmonic maps of simply connected Lorentz surfaces into general (semi-simple) semi-Riemannian symmetric spaces. (Compare to Balan and Dorfmeister [1, 3]).

For instance, in [64], the second named author gave a DPW-recipe for harmonic maps of Lorentz surfaces into the unit sphere $S^2$ and applied it to the construction of pseudospherical surfaces in Euclidean 3-space.

In [21, 22], S. Erdem studied harmonic maps of Lorentz surfaces into para-Kähler space forms, especially, paracomplex projective space $\mathbb{C}'P_n$. Since $\mathbb{C}'P_1 = S_1^2$, it seems to be interesting to construct harmonic maps of Lorentz surfaces into paracomplex projective spaces via the loop group approach.

Conformal Lorentz harmonic maps, i.e., timelike minimal immersions into semi-simple semi-Riemannian symmetric spaces have some geometric applications. For example, M. Kimura gave a construction of 3-dimensional minimal submanifolds in hyperbolic space by using timelike minimal surfaces into certain indefinite Grassmannian manifolds [37]. Burstall and Hertrich-Jeromin [8] proved Ruh–Vilms type theorems for Lie sphere geometry and projective differential geometry. They showed that appropriate Gauss maps for Lie-minimal or projective minimal surfaces are Lorentz harmonic maps into certain indefinite Grassmannian manifolds.

The following result can be deduced in much the same way as Theorem 2.5 in [18].

COROLLARY 5.4. *Let* $\psi^1, \psi^2 \colon \mathbb{D} \to S_1^2$ *be harmonic maps derived from normalized potentials* $\{\xi_1', \xi_1''\}$ *and* $\{\xi_2', \xi_2''\}$, *respectively. Then* $\psi_1$ *is congruent to* $\psi_2$ *in* $S_1^2$, *that is,* $\psi_2 = \mathrm{Ad}(\gamma)\psi_1$ *for some* $\gamma \in G$ *if and only if*

$$\xi_2' = \mathrm{Ad}(\gamma^{-1})\xi_1', \qquad \xi_2'' = \mathrm{Ad}(\gamma^{-1})\xi_1''.$$

## 6. Timelike Minimal Surfaces

### 6.1. CLASSICAL WEIERSTRASS FORMULA

Let $\varphi \colon M \to \mathbb{E}_1^3$ be a timelike surface with Gauss map $\psi$. Then the Weingarten formula implies

$$\Box \varphi = 2H\psi.$$

Thus, $M$ is minimal if and only if $\varphi$ is a vector-valued Lorentz harmonic function. Hence the immersion $\varphi$ can be written locally:

$$\varphi(u, v) = X(u) + Y(v)$$

as a sum of two curves $X(u)$ and $Y(v)$. By computing the first fundamental form of the right-hand side, we obtain

$$\langle X_u, X_u \rangle = \langle Y_v, Y_v \rangle = 0, \tag{6.1}$$

and $X_u$ and $Y_v$ are linearly independent. Hence $X(u)$ and $Y(v)$ are *null curves* in $\mathbb{E}_1^3$. For Frenet–Serret formula for null curves, we refer to Appendix C.

Here we arrive at the classical representation formula:

PROPOSITION 6.1 ([50], Theorem 3.5). *Let $\varphi\colon M \to \mathbb{E}_1^3$ be a timelike minimal surface. Then $\varphi$ is expressed locally as a sum of null curves*:

$$\varphi(u, v) = X(u) + Y(v).$$

*The velocity vector fields of the null curves $X(u)$ and $Y(v)$ are linearly independent.*

*Conversely, let $X(u)$ and $Y(v)$ be null curves defined on open intervals $I_u$ and $I_v$, respectively. Assume that the velocity vector fields $X_u$ and $Y_v$ are linearly independent. Then $\varphi(u, v) = X(u) + Y(v)$ is a timelike minimal immersion of $(I_u \times I_v, \mathrm{d}u\,\mathrm{d}v)$ into $\mathbb{E}_1^3$ with metric $\mathrm{I} = 2\langle X_u, Y_v \rangle \,\mathrm{d}u\,\mathrm{d}v$.*

Let $\varphi\colon M \to \mathbb{E}_1^3$ be a timelike surface. Then one can easily verify that the associated family of $\varphi$ is given by ([67], p. 186):

$$\varphi_\lambda(u, v) = \lambda X(u) + \lambda^{-1} Y(v), \quad \lambda \in \mathbb{R}^\times.$$

Since $\varphi_{-\lambda} = -\varphi_\lambda$, it is sufficient to investigate the associated families defined for *positive* $\lambda$. For this reason, we now restrict $\lambda$ to $\mathbb{R}^+$.

The *conjugate timelike minimal surface* $\hat\varphi$ is defined by

$$\hat\varphi(u, v) = X(u) - Y(v).$$

In the Euclidean 3-space, the conjugate minimal surface of a minimal surface belongs to the associated family of the original minimal surface. Analogously, the conjugate maximal surface of a spacelike maximal surface in Minkowski 3-space belongs to the associated family of the original one, too. However, in the timelike case, the conjugate timelike minimal surface does **not** belong to the associated family [67]. One can see that the conjugate timelike minimal surface $\hat\varphi$ is *anti isometric* to $\varphi$.

In fact, with respect to the Lorentz isothermal coordinates $(x, y)$, we have

$$\frac{\partial \hat\varphi}{\partial x} = \frac{\partial \varphi}{\partial y}, \qquad \frac{\partial \hat\varphi}{\partial y} = \frac{\partial \varphi}{\partial x}.$$

These formulas imply

$$\langle \hat\varphi_* X, \hat\varphi_* Y \rangle = -\langle \varphi_* X, \varphi_* Y \rangle$$

for all $X, Y \in \mathfrak{X}(M)$. Thus $(M, \varphi)$ is anti isometric to $(M, \hat\varphi)$.

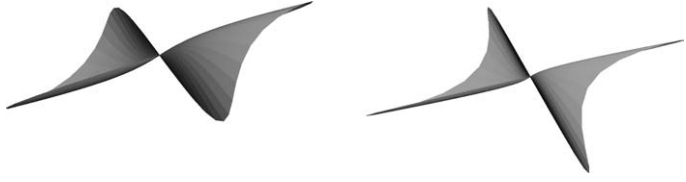

*Figure 1.* Lorentz catenoid.

COROLLARY 6.2. *Let $\varphi\colon M \to \mathbb{E}_1^3$ be a timelike minimal surface. Then, the associated family of $\varphi$ can be rewritten as*

$$\varphi_\lambda(u,v) = \cosh\theta\, \varphi(u,v) + \sinh\theta\, \hat{\varphi}(u,v), \quad \lambda = e^\theta \in \mathbb{R}^+.$$

*The tangent planes $T_{\varphi_\lambda(x,y)}M$ are mutually parallel.*

EXAMPLE 6.3 (Lorentz catenoid). Let us consider the following timelike minimal surface:

$$\varphi(u,v) = X(u) + Y(v),$$

$$X(u) = \left(\frac{\sinh\sqrt{2}u}{2}, \frac{\cosh\sqrt{2}u}{2}, \frac{u}{\sqrt{2}}\right),$$

$$Y(v) = \left(-\frac{\sinh\sqrt{2}v}{2}, -\frac{\cosh\sqrt{2}v}{2}, -\frac{v}{\sqrt{2}}\right).$$

Then

$$\varphi(u,v) = \left(\sinh\left(\frac{u-v}{\sqrt{2}}\right)\cosh\left(\frac{u+v}{\sqrt{2}}\right),\right.$$
$$\left.\sinh\left(\frac{u-v}{\sqrt{2}}\right)\sinh\left(\frac{u+v}{\sqrt{2}}\right), \frac{u-v}{\sqrt{2}}\right),$$

$$\hat{\varphi}(u,v) = \left(\cosh\left(\frac{u-v}{\sqrt{2}}\right)\sinh\left(\frac{u+v}{\sqrt{2}}\right),\right.$$
$$\left.\cosh\left(\frac{u-v}{\sqrt{2}}\right)\cosh\left(\frac{u+v}{\sqrt{2}}\right), \frac{u+v}{\sqrt{2}}\right).$$

The surface $\varphi$ is a timelike minimal surface of revolution with timelike profile curve and spacelike axis which is called the Lorentz catenoid (or timelike catenoid) with spacelike axis. The conjugate timelike minimal surface $\hat{\varphi}$ is a timelike ruled minimal surface which is called the *timelike helicoid* with spacelike axis.

Figure 1 shows the timelike catenoid viewed from different positions.

*Remark* 6.4. Timelike minimal surfaces of revolution are classified as follows (see [50], Proposition 1.21):

The nonplanar timelike minimal surfaces of revolution are congruent to one of the following:

(1) (Spacelike axis and timelike profile curve)

$$\varphi(s,t) = \left(\frac{\sinh(as+b)}{a}\cosh t, \frac{\cosh(as+b)}{a}\sinh t, s\right).$$

(2) (Spacelike axis and Euclidean profile curve)

$$\varphi(s,t) = \left(\frac{\cosh(as+b)}{a}\sinh t, \frac{\cosh(as+b)}{a}\cosh t, s\right).$$

(3) (Timelike axis)

$$\left(s, \frac{\sin(as+b)}{a}\cos t, \frac{\sin(as+b)}{a}\sin t\right).$$

(4) (Null axis)

$$(-a^2 s + b)^{1/3} L_1 + \left(-\frac{t^2}{2}(-as^2+b)^{1/3} + s\right) L_2 + t(-as^2+b)^{1/3} L_3.$$

Here $\{L_1, L_2, L_3\}$ is a null frame given by

$$L_1 = \begin{pmatrix} -1/\sqrt{2} \\ 1/\sqrt{2} \\ 0 \end{pmatrix}, \qquad L_2 = \begin{pmatrix} 1/\sqrt{2} \\ 1/\sqrt{2} \\ 0 \end{pmatrix}, \qquad L_3 = \begin{pmatrix} 0 \\ 0 \\ 1 \end{pmatrix}.$$

In the same way (see Theorem 2 in [65]), Van de Woestijine independently classified timelike minimal surfaces of revolution. Van de Woestijine calls the timelike catenoid with null axis the *surface of Enneper of the* 3rd *kind*.

Timelike minimal ruled surfaces have been classified in [65].

## 7. Timelike Minimal Surfaces via the Loop Group Method

7.1. NORMALIZED POTENTIALS FOR TIMELIKE MINIMAL SURFACES

The Lorentz conformal structure of the Lorentz symmetric space $S_1^2 = \mathrm{SL}_2\mathbb{R}/\mathrm{SO}_1(2)$ is given explicitly by

$$T_{\mathbf{k}'}S_1^2 = \mathfrak{p} = (T_{\mathbf{k}'}S_1^2)_+ \oplus (T_{\mathbf{k}'}S_1^2)_-,$$

$$(T_{\mathbf{k}'}S_1^2)_+ = \mathbb{R}\begin{pmatrix} 0 & 0 \\ 1 & 0 \end{pmatrix}, \qquad (T_{\mathbf{k}'}S_1^2)_- = \mathbb{R}\begin{pmatrix} 0 & 1 \\ 0 & 0 \end{pmatrix}.$$

Hence a map $\psi\colon \mathbb{D} \to S_1^2$ is Lorentz conformal if and only if one framing, and hence, in turn, every framing $\Psi\colon \mathbb{D} \to \mathrm{SL}_2\mathbb{R}$ has the form

$$\Psi^{-1}\mathrm{d}\Psi = \alpha'_{\mathfrak{p}} + \alpha_{\mathfrak{k}} + \alpha''_{\mathfrak{p}}$$

with

$$\alpha'_{\mathfrak{p}} = \begin{pmatrix} 0 & 0 \\ p'(u) & 0 \end{pmatrix}\mathrm{d}u, \qquad \alpha''_{\mathfrak{p}} = \begin{pmatrix} 0 & -p''(v) \\ 0 & 0 \end{pmatrix}\mathrm{d}v.$$

Here $p' = p'(u)$ and $p'' = p''(v)$ are smooth functions. The following two results can be proved in much the same way as in [18].

THEOREM 7.1. *Let $\varphi\colon \mathbb{D} \to \mathbb{E}_1^3$ be a timelike immersion and $\psi\colon \mathbb{D} \to S_1^2$ its Gauss map. Then the following are mutually equivalent*:

(1) *$\varphi$ is minimal, i.e., $H = 0$;*
(2) *$\psi$ is Lorentz-conformal;*
(3) *the normalized potential $\{\xi', \xi''\}$ has the form*

$$\xi' = \lambda \begin{pmatrix} 0 & 0 \\ p'(u) & 0 \end{pmatrix}\mathrm{d}u, \qquad \xi'' = \lambda^{-1} \begin{pmatrix} 0 & -p''(v) \\ 0 & 0 \end{pmatrix}\mathrm{d}v.$$

COROLLARY 7.2. *Let $\psi\colon \mathbb{D} \to S_1^2$ be a Lorentz conformal map and $\psi_\lambda$ be its associated family. Then*

$$\psi_\lambda = \lambda^{1/2}\psi.$$

Here we interpret $\lambda \in \mathbb{R}^+$ as $\lambda = \cosh\theta \mathbf{1} + \sinh\theta \mathbf{k}' \in K$.

### 7.2. GEOMETRIC MEANING OF THE POTENTIALS

In this subsection, we investigate differential geometric meaning of the normalized potentials of timelike minimal surfaces.

We start with the potential:

$$\xi' = \lambda \begin{pmatrix} 0 & 0 \\ p'(u) & 0 \end{pmatrix}\mathrm{d}u, \qquad \xi'' = \lambda^{-1} \begin{pmatrix} 0 & -p''(v) \\ 0 & 0 \end{pmatrix}\mathrm{d}v.$$

Then we solve the initial value problems:

$$\mathrm{d}'\Phi' = \Phi'\xi', \qquad \mathrm{d}''\Phi'' = \Phi''\xi'', \qquad \Phi'(u=0) = \Phi''(v=0) = \mathbf{1}.$$

The solution $(\Phi', \Phi'')$ is

$$\Phi'(u, \lambda) = \begin{pmatrix} 1 & 0 \\ \lambda q(u) & 1 \end{pmatrix}, \qquad \Phi''(v, \lambda^{-1}) = \begin{pmatrix} 1 & -\lambda^{-1} r(v) \\ 0 & 1 \end{pmatrix},$$

where $q(u) = \int_0^u p'(u)\,du$ and $r(v) = \int_0^v p''(v)\,dv$.

Then the Iwasawa splitting of $(\Phi', \Phi'')$ is

$$(\Phi', \Phi'') = (\Phi, \Phi)(L_-^{-1}, L_+^{-1}),$$

$$\Phi(u, v; \lambda) = \frac{1}{\sqrt{1+qr}} \begin{pmatrix} 1 & -\lambda^{-1} r(v) \\ \lambda q(u) & 1 \end{pmatrix},$$

$$L_-^{-1} = \frac{1}{\sqrt{1+q(u)r(v)}} \begin{pmatrix} 1+q(u)r(v) & \lambda^{-1} r(v) \\ 0 & 1 \end{pmatrix},$$

$$L_+^{-1} = \frac{1}{\sqrt{1+q(u)r(v)}} \begin{pmatrix} 1 & 0 \\ -\lambda q(u) & 1+q(u)r(v) \end{pmatrix}.$$

Hence the harmonic (Gauss) map is

$$N_\lambda = \mathrm{Ad}(\Phi)\mathbf{k}' = \frac{1}{1+qr} \begin{pmatrix} -\lambda q(u) + \lambda^{-1} r(v) \\ -\lambda q(u) - \lambda^{-1} r(v) \\ 1 - q(u)r(v) \end{pmatrix}.$$

The $u$-derivative of $\Phi$ is given by

$$\Phi^{-1}\Phi_u = \frac{1}{1+q(u)r(v)} \begin{pmatrix} q_u(u)r(v)/2 & 0 \\ \lambda q_u(u) & -q_u(u)r(v)/2 \end{pmatrix}. \tag{7.1}$$

Similarly, the $v$-derivative is given by

$$\Phi^{-1}\Phi_v = \frac{1}{1+q(u)r(v)} \begin{pmatrix} -q(u)r_v(v)/2 & -\lambda^{-1} r_v(v) \\ 0 & q(u)r_v(v)/2 \end{pmatrix}. \tag{7.2}$$

We determine the fundamental forms of the corresponding timelike minimal surface $\varphi \colon \mathbb{D} \to \mathbb{E}_1^3$. Take a framing $\Phi$ as in (4.10) We denote by the same letter $\Phi$ the extended framing through the framing (4.10). We call the resulting extended framing the *coordinate extended framing*. More explicitly, the coordinate extended framing is defined by

$$\mathrm{Ad}(\Phi)\left(\frac{\lambda}{2}(-\mathbf{i}+\mathbf{j}'), \frac{\lambda^{-1}}{2}(\mathbf{i}+\mathbf{j}'), \mathbf{k}'\right) = (\mathrm{e}^{-\omega/2}\varphi_u, \mathrm{e}^{-\omega/2}\varphi_v, N). \tag{7.3}$$

Direct calculations show

$$\Phi^{-1}\Phi_u = \begin{pmatrix} a & b \\ c & d \end{pmatrix} \begin{pmatrix} \omega_u/4 & 0 \\ \lambda Q \mathrm{e}^{-\omega/2} & -\omega_u/4 \end{pmatrix}, \tag{7.4}$$

$$\Phi^{-1}\Phi_v = \begin{pmatrix} a & b \\ c & d \end{pmatrix} \begin{pmatrix} -\omega_v/4 & -\lambda^{-1} R \mathrm{e}^{-\omega/2} \\ 0 & \omega_v/4 \end{pmatrix}. \tag{7.5}$$

Comparing these with (7.1)–(7.2), we get

$$\exp\{\omega(u,v)/2\} = 1 + q(u)r(v), \qquad Q(u) = q_u, \qquad R(v) = r_v.$$

## 8. The Classical Formulas, Revisited

In this section, we derive the classical representation formulas for timelike minimal surfaces from the extended frame equations. Express the coordinate extended framing $\Phi$ as $\Phi = \begin{pmatrix} a & b \\ c & d \end{pmatrix}$. Then $\Phi$ satisfies

$$\frac{\partial}{\partial u}\begin{pmatrix} a & b \\ c & d \end{pmatrix} = \begin{pmatrix} a & b \\ c & d \end{pmatrix}\begin{pmatrix} \omega_u/4 & 0 \\ \lambda Q e^{-\omega/2} & -\omega_u/4 \end{pmatrix},$$

$$\frac{\partial}{\partial v}\begin{pmatrix} a & b \\ c & d \end{pmatrix} = \begin{pmatrix} a & b \\ c & d \end{pmatrix}\begin{pmatrix} -\omega_v/4 & -\lambda^{-1} R e^{-\omega/2} \\ 0 & \omega_v/4 \end{pmatrix}.$$

From these equations, we deduce

$$a = e^{-\omega/4} s_1(u), \qquad b = e^{-\omega/4} t_1(v),$$
$$c = e^{-\omega/4} s_2(u), \qquad d = e^{-\omega/4} t_2(v).$$

Inserting these into (7.3), we obtain

$$\varphi_u = \left(\frac{-1}{2}(s_1(u)^2 + s_2(u)^2), \frac{1}{2}(s_1(u)^2 - s_2(u)^2), s_1(u)s_2(u)\right),$$

$$\varphi_v = \left(\frac{1}{2}(t_1(v)^2 + t_2(v)^2), \frac{-1}{2}(t_1(v)^2 - t_2(v)^2), -t_1(v)t_2(v)\right).$$

Thus $\varphi$ is congruent to

$$\varphi = \int_0^u \left(\frac{-1}{2}(s_1(u)^2 + s_2(u)^2), \frac{1}{2}(s_1(u)^2 - s_2(u)^2), s_1(u)s_2(u)\right) du +$$
$$+ \int_0^v \left(\frac{1}{2}(t_1(v)^2 + t_2(v)^2), \frac{-1}{2}(t_1(v)^2 - t_2(v)^2), -t_1(v)t_2(v)\right) dv$$

up to translations. The Gauss map is

$$\psi = e^{-\omega/2}\begin{pmatrix} -s_1(u)t_1(v) - s_2(u)t_2(v) \\ s_1(u)t_1(v) - s_2(u)t_2(v) \\ s_1(u)t_2(v) + s_2(u)t_1(v) \end{pmatrix}.$$

The projected Gauss map $\mathcal{P}_+ \circ \psi$ is given by

$$\frac{1}{2}\left(-\frac{t_1(v)}{t_2(v)} - \frac{s_2(u)}{s_1(u)}, \frac{t_1(v)}{t_2(v)} + \frac{s_2(u)}{s_1(u)}\right).$$

Comparing and evaluating the following equation at $\lambda = 1$:

$$\Phi = \frac{1}{\sqrt{1+qr}} \begin{pmatrix} 1 & -\lambda^{-1}r \\ \lambda q & 1 \end{pmatrix} = e^{-\omega/4} \begin{pmatrix} s_1 & t_1 \\ s_2 & t_2 \end{pmatrix},$$

we get

$$s_1(u) = 1, \qquad s_2(u) = q(u), \qquad t_1(v) = -r(v), \qquad t_2(v) = 1.$$

Thus we have

$$\varphi_u = \left(-\frac{1}{2}(1+q^2), \frac{1}{2}(1-q^2), q\right), \qquad \varphi_v = \left(\frac{1}{2}(1+r^2), \frac{1}{2}(1-r^2), r\right).$$

This is the *normalized Weierstrass formula* obtained by Magid. See Theorem 4.3 and p. 456, Notes 2 in [49].

In [49], geometric meaning of the integrand data $q$ and $r$ are not clarified. Our discussion shows that the data $q$ and $r$ are primitive functions of the coefficients of the Hopf differential.

Finally, we investigate the normalized Weierstrass formula by using the theory of null Frenet curves.

The null curves

$$X(u) = \int_0^u \left(-\frac{1}{2}(1+q(u)^2), \frac{1}{2}(1-q(u)^2), q(u)\right) du,$$

$$Y(v) = \int_0^v \left(\frac{1}{2}(1+r(v)^2), \frac{1}{2}(1-r(v)^2), r(v)\right) dv$$

have null Frenet frame fileds

$$A^X(u) = \left(-\frac{1}{2}(1+q(u)^2), \frac{1}{2}(1-q(u)^2), q(u)\right),$$

$$B^X(u) = (1, 1, 0), \qquad C^X(u) = (-q(u), -q(u), 1),$$

$$A^Y(v) = \left(\frac{1}{2}(1+r(v)^2), \frac{1}{2}(1-r(v)^2), r(v)\right),$$

$$B^Y(v) = (-1, 1, 0), \qquad C^Y(v) = (r(v), -r(v), 1),$$

respectively. The curvature functions of $X$ and $Y$ are

$$k_1^X(u) = q_u = Q, \qquad k_1^Y(v) = r_v = R, \qquad k_2^X(u) = k_2^Y(v) = 0.$$

*Remark 8.1.* A timelike surface in $\mathbb{E}_1^3$ is said to be *isothermic* if it is covered by isothermic coordinate system [26]. A null coordinate system $(u, v)$ is said to be an isothermic coordinate system if $Q = R$.

For every simply connected timelike isothermic surface $\varphi \colon \mathbb{D} \to \mathbb{E}_1^3$, there exists a timelike isothermic surface $\varphi^* \colon \mathbb{D} \to \mathbb{E}_1^3$ with fundamental forms:

$$\mathrm{I}^* = \mathrm{e}^{-\omega}\,\mathrm{d}u\,\mathrm{d}v, \qquad Q^* = R^* = H/2, \qquad H^* = 2Q.$$

The timelike isothermic surface $\varphi^*$ is called the *dual surface* of $\varphi$. If $(\mathbb{D}, \varphi)$ is a timelike isothermic minimal surface, then its dual surface is the Gauss map $\psi \colon \mathbb{D} \to S_1^2$. Conversely, let $\psi \colon \mathbb{D} \to S_1^2$ be a Lorentz-conformal map. Then $\psi$ is isothermic and its dual surface is minimal. Moreover, $\psi^*$ is given explicitly by the normalized Weierstrass formula discussed above.

## 9. Examples

We conclude this paper with exhibiting fundamental examples.

EXAMPLE 9.1. We start with a normalized potential given by $p' = \varepsilon = \pm 1$, $p'' = 1$. The projected Gauss map determined by the potential is

$$-\frac{\varepsilon}{2}(u - \varepsilon v, u + \varepsilon v). \tag{9.1}$$

Now we give a timelike minimal surface with projected Gauss map (9.1). By the classical formula:

$$\varphi_u^{(\varepsilon)} = \left(-\frac{1}{2}(1+q^2), \frac{1}{2}(1-q^2), q\right),$$
$$\varphi_v^{(\varepsilon)} = \left(\frac{1}{2}(1+r^2), \frac{1}{2}(1-r^2), r\right),$$

we obtain the following immersion:

$$\varphi^{(\varepsilon)}(u, v) = X(u) + Y(v),$$

where

$$X(u) = \frac{1}{2}\left(-\left(u + \frac{u^3}{3}\right), u - \frac{u^3}{3}, \varepsilon u^2\right),$$

$$Y(v) = \frac{1}{2}\left(v + \frac{v^3}{3}, v - \frac{v^3}{3}, v^2\right).$$

These formulas show that the parameters $u$ and $v$ are pseudoarc parameters of $X(u)$ and $Y(v)$, respectively. And hence $X(u)$ and $Y(v)$ are null helices of curvature 0.

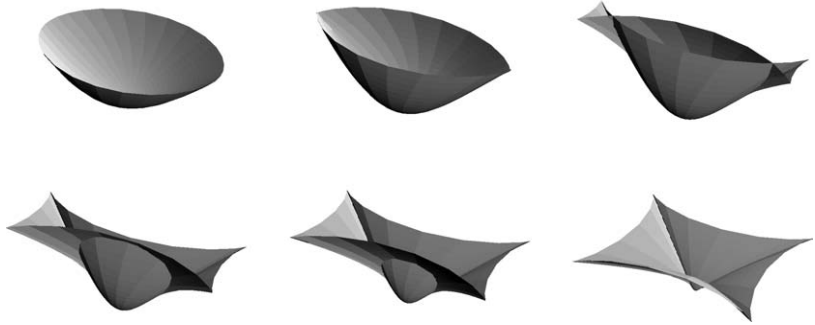

*Figure 2.* Lorentz cousin of Enneper's surface (isothermic type).

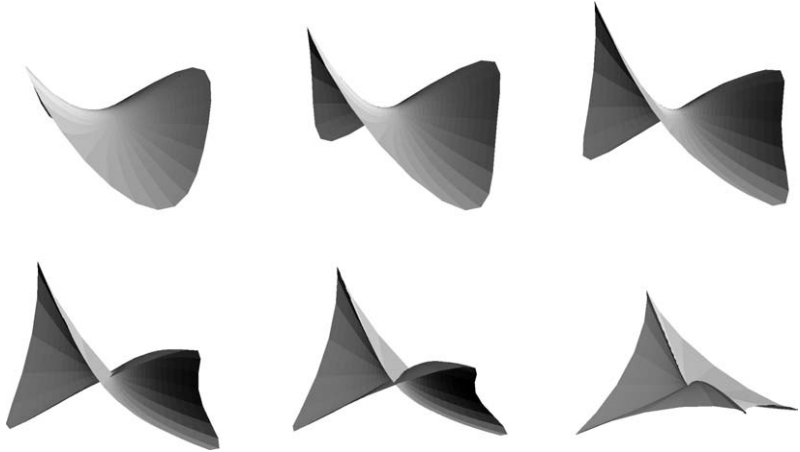

*Figure 3.* Lorentz cousin of Enneper's surface (anti-isothermic type).

The timelike surface $\varphi^{(1)}$ may be regarded as "Lorentzian cousin" of Enneper's minimal surface. The Hopf differentials and the Gaussian curvature of $\varphi^{(\varepsilon)}$ are

$$Q = \varepsilon, \qquad R = 1, K = -4\varepsilon(1 + \varepsilon u v)^{-4}.$$

The metric of $\varphi^{(\varepsilon)}$ is

$$\mathrm{I} = (1 + \varepsilon u v)^2 \, \mathrm{d}u \, \mathrm{d}v.$$

Thus $\varphi^{(1)}$ is isothermic, i.e., it has real distinct principal curvatures. On the contrary, $\varphi^{(-1)}$ has no Euclidean counterpart, since it is anti isothermic. Namely $\varphi^{(-1)}$ has imaginary principal curvatures. Both the surfaces are foliated by *null helices* of curvature 0. The Enneper's minimal surface in Euclidean 3-space does not have such a property.

EXAMPLE 9.2. Next we consider a normalized potential with $p' \neq 0$ and $p'' = 0$. By the classical formula, we obtain

$$\varphi(u, v) = X(u) + Y(v),$$

$$X(u) = \int_0^u \left(-\frac{1}{2}(1+q(u)^2), \frac{1}{2}(1-q(u)^2), q(u)\right) du,$$

$$Y(v) = v\left(\frac{1}{2}, \frac{1}{2}, 0\right).$$

This explicit representation shows that $\varphi$ is a ruled surface. The null Frenet frame filed $(A^X, B^X, C^X)$ along the null curve $X(u)$ is given by

$$A^X(u) = X_u(u) = \left(-\frac{1}{2}(1+q(u)^2), \frac{1}{2}(1-q(u)^2), q(u)\right),$$

$$B^X(u) = (1, 1, 0), \qquad C^X(u) = (-q, -q, 1).$$

Hence we notice the $Y(v) = vB^X(u)/2$. Thus $\varphi$ is the $B$-scroll of $X(u)$ (see Appendix C). Note that the unit normal vector field of $\varphi$ is $N = C$. The second fundamental form of $\varphi$ is described as $Q = q_u(u)$, $R = 0$, $K = H = 0$. Hence $\varphi$ has principal curvature 0 with multiplicity 2. However $\varphi$ is not totally geodesic since $q_u \neq 0$ by our assumption $p' \neq 0$. This timelike surface has no Euclidean counterpart. Reparametrize $(u, v)$ so that $Q = 1$. Without loss of generality, we may assume that $q(u) = u$. Then $u$ is the pseudoarc parameter of $X(u)$. Then $\varphi$ is parametrized as

$$\varphi(u, v) = \frac{1}{2}\left(-\frac{u^3}{3} - u + v, -\frac{u^3}{3} + u + v, u^2\right).$$

Hence $\varphi$ is a cylinder over a parabora $2\xi_3 = (-\xi_1 + \xi_2)^2$. P. Mira and J. A. Pastor called a minimal $B$-scroll a *parabolic null cylinder* [54]. In [54], Mira and Pastor obtained some characterizations of the parabolic null cylinder.

## 10. Concluding Remarks

*Remark 10.1.* The Gauss–Codazzi equations (4.6)–(4.7)–(4.8) imply the following fact [25, 26]:

PROPOSITION 10.2. *Let $\varphi \colon \mathbb{D} \to \mathbb{E}_1^3$ be a timelike immersion of a simply connected Lorentz surface with global null coordinate system with constant mean curvature $H_0$. Then there exists a timelike immersions*

$$\varphi^S \colon \mathbb{D} \to S_1^3, \qquad \varphi^H \colon \mathbb{D} \to H_1^3$$

*of $\mathbb{D}$ into the de Sitter 3-space $S_1^3$ and the anti de Sitter 3-space $H_1^3$ such that both $\varphi^S$ and $\varphi^H$ are isometric to $\varphi$ and have constant mean curvatures $H_S$ and $H_H$ which satisfy the following relations:*

$$H_0^2 = H_S^2 + 1 = H_H^2 - 1.$$

*The correspondences between $\varphi$, $\varphi^S$ and $\varphi^H$ are called the Lawson–Guichard correspondences.*

In particular, there exists (local) isometric correspondences between timelike minimal surfaces in $\mathbb{E}^3_1$ and timelike surfaces of constant mean curvature $\pm 1$ in $H^3_1$. This correspondene induces a Weierstrass-type representation for timelike surfaces with mean curvature $\pm 1$ in $H^3_1$ (cf. [33]).

*Remark 10.3.* Timelike minimal surfaces [resp. spacelike maximal surfaces] in $H^3_1$ are conformal harmonic maps of Lorentz [resp. Riemann] surfaces into $H^3_1$. By using $\mathrm{SL}_2\mathbb{R}$-model of $H^3_1$, we regard timelike minimal [resp. spacelike maximal] surfaces as "harmonic maps" into the noncompact simple Lie group $\mathrm{SL}_2\mathbb{R}$. We hope a general scheme by Balan and Dorfmeister [1] would enable us to get substantial progress in the study of such surfaces in $H^3_1$. Note that study of timelike minimal or spacelike maximal surfaces in $H^3_1$ has some motivations from particle physics [4, 5].

*Remark 10.4.* Clelland [10] investigated Bäcklund-type transformations (line congruences) for timelike constant mean curvature surfaces. Her transformations can be described as "dressing actions" of loop groups. This topic will be discussed in a future work.

**Acknowledgements**

The authors would like to thank Professor Josef Dorfmeister for fruitful discussions and constant encouragement. In particular, we are most indebted to him for helping us prove the Birkhoff decomposition theorem, which is presented in the Appendix. The authors would also like to thank Shoichi Fujimori for contributing the graphics in this paper.

**Appendix A. The Birkhoff Decomposition Theorem**

In this appendix, we give a proof of the Birkhoff decomposition theorem (Theorem 5.1) for $\tilde{\Lambda} G$ for completeness.

Let us denote the polynomial loop algebra of $\mathfrak{g} = \mathfrak{sl}_2\mathbb{R}$ by $\Lambda_{\mathrm{pol}}\mathfrak{g}$:

$$\Lambda_{\mathrm{pol}}\mathfrak{g} = \left\{ \xi(\lambda) = \sum_{\mathrm{finite}} \xi_j \lambda^j : S^1 \to \mathfrak{g} \right\}.$$

We introduce the norm $|\cdot|$ for $\mathfrak{g}$:

$$|A| = \max_{j=1,2}\left\{ \sum_{i=1}^{2} |a_{ij}|,\ A = (a_{ij}) \in \mathfrak{g} \right\}.$$

We extend this norm to $\Lambda_{\mathrm{pol}}\mathfrak{g}$ in the following way:

$$\|\xi\| = \sum |\xi_j|, \quad \xi = \sum \xi_j \lambda^j.$$

Denote the Banach-completion of $\Lambda_{\mathrm{pol}}\mathfrak{g}$ by $\Lambda\mathfrak{g}$. The resulting Banach Lie algebra is called the (free) *loop algebra* of $\mathfrak{g}$.

One can analogously define a connected Banach Lie group $\Lambda G$ whose Lie algebra is $\Lambda\mathfrak{g}$.

Moreover loop groups for the complexified group $G^{\mathbb{C}} = \mathrm{SL}_2\mathbb{C}$ are similarly defined.

Here we recall the classical Birkhoff Decomposition Theorem for the free loop group $\Lambda G^{\mathbb{C}}$.

THEOREM A (Birkhoff decomposition of $\Lambda G^{\mathbb{C}}$ [58]).

$$\Lambda G^{\mathbb{C}} = \bigsqcup_{w \in \mathcal{T}} \Lambda^- G^{\mathbb{C}} \cdot w \cdot \Lambda^+ G^{\mathbb{C}},$$

*Here $\mathcal{T}$ denotes the group of homomorphisms from $S^1$ into the subgroup of diagonal matrices of $\mathrm{SL}_2\mathbb{C}$, that is,*

$$\mathcal{T} = \left\{ \begin{pmatrix} \lambda^a & 0 \\ 0 & \lambda^{-a} \end{pmatrix} \,\middle|\, a > 0 \right\}.$$

*Moreover, the multiplication maps*

$$\Lambda_*^- G^{\mathbb{C}} \times \Lambda^+ G^{\mathbb{C}} \to \Lambda G^{\mathbb{C}}, \qquad \Lambda_*^+ G^{\mathbb{C}} \times \Lambda^- G^{\mathbb{C}} \to \Lambda G^{\mathbb{C}}$$

*are diffeomorphisms onto the open dense subsets $\mathcal{B}_\Lambda^\circ(-,+)$ and $\mathcal{B}_\Lambda^\circ(+,-)$ of $\Lambda G^{\mathbb{C}}$, called the big cells of $\Lambda G^{\mathbb{C}}$. In particular if $\gamma$ is an element of $\mathcal{B}_\Lambda = \mathcal{B}_\Lambda^\circ(-,+) \cap \mathcal{B}_\Lambda^\circ(+,-)$, then $\gamma$ has unique decompositions*:

$$\gamma = \gamma_- \cdot \ell_+ = \gamma_+ \cdot \ell_-, \quad \gamma_\pm \in \Lambda_*^\pm G^{\mathbb{C}}, \ell_\pm \in \Lambda^\pm G^{\mathbb{C}}.$$

*Here the subgroups $\Lambda_*^\pm G^{\mathbb{C}}$ are defined by*

$$\Lambda_*^- G^{\mathbb{C}} = \left\{ \gamma \in \Lambda^- G^{\mathbb{C}} \mid \gamma(\lambda) = \mathbf{1} + \sum_{k \leqslant -1} \gamma_k \lambda^k \right\},$$

$$\Lambda_*^+ G^{\mathbb{C}} = \left\{ \gamma \in \Lambda^+ G^{\mathbb{C}} \mid \gamma(\lambda) = \mathbf{1} + \sum_{k \geqslant 1} \gamma_k \lambda^k \right\}.$$

Let $\tilde{\Lambda} G$ be the subset of $\Lambda G$ whose elements admit analytic continuations to $\mathbb{C}^\times$ as in Section 2. For $\tilde{\Lambda} G$ we will use the topology induced from $\Lambda G$. Then we obtain Theorem 5.1.

*Proof of the Birkhoff Decomposition Theorem 5.1.* Let $g \in \tilde{\Lambda} G$ with expansion $g(\lambda) = \sum g_j \lambda^j$.

Note that the coefficients $g_j$ in the expansion of $g$ are real.

Over the unit circle $S^1$, we obtain a decomposition $g = g_- \cdot w \cdot g_+$ by the classical Birkhoff Decomposition Theorem. It remains to show that actually every factor of $g$ defines an element in $\tilde{\Lambda} G$.

First we show all factors of $g$ are in $\Lambda G$. To see this we introduce the automorphism $\kappa$ of $\Lambda G^{\mathbb{C}}$ which is defined by

$$\kappa(\gamma)(\lambda) := \sum \bar{\gamma}_j \lambda^j$$

for every

$$\gamma(\lambda) = \sum \gamma_j \lambda^j \in \Lambda G^{\mathbb{C}}.$$

It is obvious to verify that $\kappa$ leaves $\Lambda^+ G^{\mathbb{C}}$ and $\Lambda^- G^{\mathbb{C}}$ invariant and fixes all $w \in \mathcal{T}$. Moreover, $\Lambda G$ is the fixed point set of $\kappa$. Thus $g = \kappa(g) = \kappa(g_-) \cdot w \cdot \kappa(g_+)$. The classical Birkhoff Decomposition Theorem now implies $\kappa(g_-) = g_- \cdot v_-$ and $\kappa(g_+) = v_+ \cdot g_+$, where $v_\pm \in \Lambda^\pm G^{\mathbb{C}}$. Moreover, we have $v_- \cdot w = w \cdot v_+^{-1}$.

Applying $\kappa$ again to $\kappa(g_-) = g_- v_-$ and taking into account that $\kappa$ is an automorphism of order two, it follows that

$$g_- = \kappa(\kappa(g_-)) = \kappa(g_- v_-) = g_- \cdot v_- \cdot \kappa(v_-).$$

Thus we obtain $\kappa(v_-) = (v_-)^{-1}$.

Analogusly, from $\kappa(g_+) = v_+ g_-$, we obtain $\kappa(v_+) = (v_+)^{-1}$.

If $w = \mathbf{1}$, then we are done. Assume now $w \neq \mathbf{1}$. Then $v_-$ and $v_+$ are lower triangular. Moreover, the equation relating $v_-$ and $v_+$ via $w$ shows that the diagonal part of $v_-$ is independent of $\lambda$. Thus $\kappa(v_-) = (v_-)^{-1}$ shows that the diagonal entries of $v_-$ have modulus 1. Let $d$ denote the diagonal part of $v_-$ and let $\sqrt{d}$ denote its square root. Then we set $g'_- := g_- \cdot (\sqrt{d})^{-1}$ and obtain $\kappa(g'_-) = g'_- \cdot v'_-$, where $v'_-$ is lower triangualr with 1's on the diagonal. Writing now $g = g'_- \cdot w \cdot g'_+$, where $g'_+$ is defined by $g'_+ = \sqrt{d} \cdot g_+$, then we see that the corresponding $v'_+$ has 1's on the diagonal. Now we take one of the square root of $v'_-$ and multiply $g'_-$ by its inverse on the right, obtaining $g''_-$. A straightforward computation shows that $\kappa$ fixes $g''_-$, whence $g''_-$ is in $\Lambda G$. As a consequence, the corresponding $g''_+$ is also in $\Lambda G$.

Next we need to show that $g_-$ and $g_+$ are in the "$\sim$"-group.

We know from the definition that $g_-$ has a holomorphic extension to the exterior of the unit disk and is finite at $\infty$. Thus $(g_-)^{-1} \cdot g = w \cdot g_+$ has a holomorphic extension to the exterior of the unit disk. Hence $g_+$ also has a holomorphic extension to the exterior of the unit disk. Altogether, $g_+$ has a holomorphic extension to $\mathbb{C}^\times$. This proves the first part of the theorem.

For the second part we note first that the big cells are indeed dense, since the cosets involving $w \neq \mathbf{1}$ are of nonzero codimension in $\Lambda G$ (cf. [2, 16]).

The rest of the claim follows the fact the multiplication maps are induced by the diffeomorphisms $\Lambda^-_* G^{\mathbb{C}} \times \Lambda^+ G^{\mathbb{C}} \to \Lambda G^{\mathbb{C}}$ and $\Lambda^+_* G^{\mathbb{C}} \times \Lambda^- G^{\mathbb{C}} \to \Lambda G^{\mathbb{C}}$ via the loop group correspondences we presented in Section 5. □

**Appendix B. Dualities between Harmonic Maps**

In this appendix, we give a "self-duality" of harmonic maps into $S_1^2$.

In [8], Burstall and Hertrich-Jeromin pointed out a duality between Lorentz harmonic maps into symmetric spaces.

Let $\mathcal{N} = G/K$ be a semi simple semi-Riemannian symmetric space with symmetric Lie algebra decomposition $\mathfrak{g} = \mathfrak{k} \oplus \mathfrak{m}$. Then $\hat{\mathfrak{g}} := \mathfrak{k} \oplus \sqrt{-1}\mathfrak{m}$ is a symmetric Lie algebra decomposition of a second symmetric space $\hat{\mathcal{N}} = \hat{G}/\hat{K}$.

Now let $\psi \colon \mathbb{D} \to G/K$ be a harmonic map of a simply connected Lorentz surface with an extended framing $\Psi_\lambda \colon \mathbb{D} \times \mathbb{R}^+ \to G$. Since $\Psi_\lambda$ is analytic in $\lambda \in \mathbb{R}^+$, it admits an analytic continuation to $\mathbb{C}^\times$. Then for $\lambda \in \sqrt{-1}\mathbb{R}^\times$, $\Psi_\lambda$ takes value in $\hat{G}$ (by adjusting left translation by constant matrix in $\hat{G}$, if necessary).

PROPOSITION B ([8]). *Let $\mathbb{D}$ be a simply connected Lorentz surface. Then there exists a bijective correspondence between harmonic maps $\psi \colon \mathbb{D} \to \mathcal{N}$ and $\hat{\psi} \colon \mathbb{D} \to \hat{\mathcal{N}}$ modulo isometries.*

Now we specialize this duality for harmonic maps into $S_1^2$.

Let $H_{\mathbf{j}'}$ be the isotropy subgroup of $G$ at $\mathbf{j}' \in S_1^2$. Denote by $\mathfrak{h}_{\mathbf{j}'}$ the Lie algebra of $H_{\mathbf{j}'}$. The tangent space of $S_1^2$ at $\mathbf{j}'$ is denoted by $\mathfrak{m}_{\mathbf{j}'}$. Then the new Lie algebra $\hat{\mathfrak{g}} := \mathfrak{h}_{\mathbf{j}'} \oplus \sqrt{-1}\mathfrak{m}_{\mathbf{j}'}$ is

$$\left\{ \begin{pmatrix} -\sqrt{-1}\xi_3 & \xi_2 - \sqrt{-1}\xi_1 \\ \xi_2 + \sqrt{-1}\xi_1 & \sqrt{-1}\xi_3 \end{pmatrix} \right\}.$$

Namely $\hat{\mathfrak{g}}$ is the Lie algebra $\mathfrak{su}_1(2)$ of the indefinite unitary group:

$$\mathrm{SU}_1(2) = \left\{ \begin{pmatrix} \beta_1 & \bar{\beta}_2 \\ \beta_2 & \bar{\beta}_1 \end{pmatrix} \mid |\beta_1|^2 - |\beta_2| = 1,\ \beta_1, \beta_2 \in \mathbb{C} \right\}.$$

It is well known that $\mathfrak{su}_1(2)$ is isomorphic to $\mathfrak{g}$ as a real Lie algebra. The second symmetric space is $S_1^2 = \mathrm{SU}_1(2)/\mathbb{R}^\times$. Thus every harmonic map $\psi \colon \mathbb{D} \to S_1^2 = G/H_{\mathbf{j}'}$ is dual to a harmonic map into $S_1^2$.

**Appendix C. Null Curves**

Let $\gamma = \gamma(s)$ be smooth regular curve in $\mathbb{E}_1^3$. Then $\gamma$ is said to be a *null Frenet curve* if

(1) the tangent vector field of $\gamma$ is null, i.e., $\langle \gamma', \gamma' \rangle = 0$;
(2) there exists vector fields $A$, $B$ and $C$ along $\gamma$ such that

$$A = \gamma', \qquad \langle A, B \rangle = 1, \qquad \langle B, B \rangle = 0, \qquad \langle C, C \rangle = 1,$$

$$\frac{\mathrm{d}}{\mathrm{d}s}(A, B, C) = (A, B, C) \begin{pmatrix} 0 & 0 & -k_2 \\ 0 & & -k_1 \\ k_1 & k_2 & 0 \end{pmatrix}.$$

The frame filed $\mathcal{L} = (A, B, C)$ is called the *null Frenet frame field* of $\gamma$. Note that $\langle A', A' \rangle = \kappa^2 \geqslant 0$.

For a null Frenet curve $\gamma(s)$, the ruled surface

$$\varphi(s, t) := \gamma(s) + tB(s)$$

is a timelike surface with two equal real principal curvatures. In fact, the mean curvature and the Gaussian curvatures are given by $H = k_2$ and $K = k_2^2$. This ruled surface is called the *B-scroll* of $\gamma$ [29].

When $\langle A', A' \rangle > 0$, there exists a reparametrization $s = s(\sigma)$ such that $\langle \gamma_{\sigma\sigma}, \gamma_{\sigma\sigma} \rangle = 1$. Under this new parametrization, the Frenet–Serret formula becomes

$$\frac{d}{d\sigma}(A, B, C) = (A, B, C) \begin{pmatrix} 0 & 0 & -k \\ 0 & & -1 \\ 1 & k & 0 \end{pmatrix}.$$

Ferrández, Giménez and Lucas [23] called this parameter $\sigma$ and the function $k$ the *pseudoarc parameter* and *null Cartan curvature*, respectively. Thus a null Frenet curves are completely detrmined only by the null Cartan curvature $k(\sigma)$ up to Lorentz transformations.

A null Frenet curve parametrized by pseudoarc parameter is called a *null helix* if its null Cartan curvature is constant [23].

One can check that every null helix of null Cartan curvature 0 is congruent to (cf. [23]):

$$\left( -\frac{1}{2}\left( \sigma + \frac{\sigma^3}{3} \right), \frac{1}{2}\left( \sigma - \frac{\sigma^3}{3} \right), \frac{\sigma^2}{2} \right).$$